\documentclass{amsart}
\usepackage{amssymb}
\usepackage[all]{xy}

      \newenvironment{changemargin}[2]{\begin{list}{}{
         \setlength{\topsep}{0pt}\setlength{\leftmargin}{0pt}
         \setlength{\rightmargin}{0pt}
         \setlength{\listparindent}{\parindent}
         \setlength{\itemindent}{\parindent}
         \setlength{\parsep}{0pt plus 1pt}
         \addtolength{\leftmargin}{#1}\addtolength{\rightmargin}{#2}
         }\item }{\end{list}}

\newcommand{\Bc}[9]{\bibitem{#1} {#2}, \emph{#3}, in: \textbf{#4} (#5), #6 #7, #8--#9.}

\newcommand{\plim}{{p\mbox{-}\lim}}
\renewcommand{\split}{\mathsf{Split}}
\newcommand{\mx}[1]{\begin{matrix}#1\end{matrix}}
\newcommand{\gp}{\mbox{-\textit{\tiny gp}}}
\newcommand{\gpbl}{{\mbox{\textit{\tiny gp}}}}
\newcommand{\Cite}[2]{{\cite[#1]{#2}}}
\newcommand{\Ref}[1]{\textbf{[[[#1]]]}}
\newcommand{\CH}{the Continuum Hypothesis}
\newcommand{\MA}{Martin's Axiom}
\newcommand{\BL}{\B_\Lambda}

\newcommand{\CO}{C_\Omega}
\newcommand{\arx}[1]{\texttt{http://arxiv.org/abs/#1}}
\long\def\forget#1\forgotten{}

\newcommand{\cU}{\mathcal{U}}
\newcommand{\cV}{\mathcal{V}}
\newcommand{\cW}{\mathcal{W}}
\newcommand{\Increasing}{\mathcal{Z}}
\newcommand{\cP}{\mathcal{P}}

\newcommand{\fj}{\mathfrak{j}}

\newcommand{\fU}{\mathfrak{U}}
\newcommand{\fV}{\mathfrak{V}}
\newcommand{\seq}[1]{\{#1\}_{n\in\N}}

\newcommand{\B}{\mathcal{B}}

\newcommand{\BG}{\B_\Gamma}
\newcommand{\BT}{\B_\mathrm{T}}
\newcommand{\BO}{\B_\Omega}
\newcommand{\DO}{\cD_\Omega}
\newcommand{\KO}{\cK_\Omega}

\newcommand{\cD}{\mathcal{D}}
\newcommand{\cK}{\mathcal{K}}
\newcommand{\cN}{\mathcal{N}}

\newcommand{\Tau}{\mathrm{T}}

\newcommand{\cF}{\mathcal{F}}

\newcommand{\M}{\mathcal{M}}
\newcommand{\N}{\mathbb{N}}
\newcommand{\compactN}{\mathbb{N}\cup\{\infty\}}

\newcommand{\NN}{{{}^\N\N}}
\newcommand{\NcompactN}{{{}^\N(\compactN)}}

\newcommand{\roth}{P_\oo(\N)}

\renewcommand{\O}{\mathcal{O}}

\newcommand{\R}{\mathbb{R}}

\renewcommand{\u}{\mathfrak{u}}
\newcommand{\Union}{\bigcup}
\newcommand{\Z}{\mathbb{Z}}

\renewcommand{\b}{{\mathfrak b}}

\renewcommand{\c}{{\mathfrak c}}
\renewcommand{\d}{{\mathfrak d}}

\renewcommand{\i}{\item}

\newcommand{\oo}{\infty}
\newcommand{\p}{{\mathfrak p}}
\newcommand{\s}{\mathfrak{s}}
\newcommand{\w}{\omega}
\newcommand{\x}{\times}

\newcommand{\Iff}{\Leftrightarrow}

\newcommand{\nin}{\not\in}

\newcommand{\sbst}{\subseteq}

\newcommand{\sm}{\setminus}

\renewcommand{\(}{\left(}
\renewcommand{\)}{\right)}

\newcommand{\cI}{\mathcal{I}}
\newcommand{\cJ}{\mathcal{J}}

\newcommand{\cov}{\mathsf{cov}}
\newcommand{\wQN}{\mathsf{wQN}}
\newcommand{\add}{\mathsf{add}}
\newcommand{\cof}{\mathsf{cof}}
\newcommand{\cf}{\mathsf{cf}}
\newcommand{\non}{\mathsf{non}}

\newcommand{\impl}{\to}
\renewcommand{\t}{\mathfrak{t}}
\newcommand{\h}{\mathfrak{h}}

\newtheorem{thm}{Theorem}[section]

\newtheorem{prob}[thm]{Problem}
\newtheorem{proj}[thm]{Project}

\newtheorem{conj}[thm]{Conjecture}

\theoremstyle{definition}

\theoremstyle{remark}

\newcommand{\be}{\begin{enumerate}}
\newcommand{\ee}{\end{enumerate}}
\newcommand{\bi}{\begin{itemize}}
\newcommand{\ei}{\end{itemize}}

\newcommand{\sone}{\mathsf{S}_1}    \newcommand{\sfin}{\mathsf{S}_{fin}}
\newcommand{\ufin}{\mathsf{U}_{fin}}

\newcommand{\gone}{\mathsf{G}_1}    \newcommand{\gfin}{\mathsf{G}_{fin}}

\newcommand{\cl}{\overline}

\newcommand{\lft}[2]{\mathopen\ifcase#1{}\oo\or
                        \big#2\or\Big#2\else\oo\fi}
\newcommand{\ed}{\end{document}}
\newcommand{\rgt}[2]{\mathclose\ifcase#1{}\oo\or
                        \big#2\or\Big#2\else\oo\fi}

\title[SPM problems milestone]{Selection principles in mathematics:\\ A milestone of open problems}
\author{Boaz Tsaban}
\thanks{Partially supported by the Golda Meir Fund and
the Edmund Landau Center for Research in Mathematical Analysis and Related Areas,
sponsored by the Minerva Foundation (Germany).}
\address{Einstein Institute of Mathematics, Hebrew University of Jerusalem,
Givat Ram, Jerusalem 91904, Israel}
\email{tsaban@math.huji.ac.il}
\urladdr{http://www.cs.biu.ac.il/\~{}tsaban}

\keywords{%
Gerlits-Nagy $\gamma$ property,
Galvin-Miller strong $\gamma$ property,
Menger property,
Hurewicz property,
Rothberger property,
Gerlits-Nagy $(*)$ property,
Arkhangel'ski\v{i} property,
Sakai property,
Selection principles,
Infinite game theory,
Ramsey Theory.
}

\subjclass{%
Primary: 37F20; 
Secondary 26A03, 
03E75 
}

\begin{document}

\begin{abstract}
We survey some of the major open problems involving selection
principles, diagonalizations, and covering properties in topology
and infinite combinatorics. Background details, definitions and
motivations are also provided.
\\[0.1cm]
\emph{Remark.}
The paper is not updated anymore. See\\
\centerline{\texttt{http://arxiv.org/math/0609601}}
for a more up-to-date survey.
\end{abstract}
\maketitle

\begin{changemargin}{8cm}{0cm}
\noindent\emph{A mathematical discipline is a\-li\-ve and well if it has many exciting open
problems of different levels of difficulty.}

\medskip

\hfill Vitali Bergelson \cite{Bergelson}
\end{changemargin}

\tableofcontents

\section{Introduction}
The general field in which the problems surveyed below arise is traditionally called
\emph{Selection Principles in Mathematics (SPM)}.\footnote{Some other popular names are:
Topological diagonalizations, infinite-combinatorial topology.
}
It is at least as old as Cantor's works on the diagonalization argument. However,
we concentrate on the study of diagonalizations of covers of topological spaces
(and their relations to infinite combinatorics, Ramsey theory, infinite game theory,
and function spaces) since these are the parts of this quickly-growing field with which
we are more familiar. Example for an important area which is not covered here is that of
topological groups.
For problems on these and many more problems in the areas we
do consider, the reader is referred to the papers cited in Scheepers' survey paper
\cite{LecceSurvey}.%
\footnote{Other topics which are very new, such as star selection principles and uniform selection principles,
are not covered here and the reader is referred to
the other papers in \emph{this} issue of \textbf{Note di Matematica} for some information on these.}

Many mathematicians have worked in the past on specific instances of these
``topological diagonalizations'', but it was only in 1996 that Marion Scheepers'
paper \cite{coc1} established a unified framework to study all of these sorts of
diagonalizations.
This pioneering work was soon followed by a stream of papers, which seems to get
stronger with time. A new field was born.

As always in mathematics, this systematic approach made it possible to
generalize and understand to a much deeper extent existing results,
which were classically proved using ad-hoc methods and ingenious arguments
which were re-invented for each specific question.

However, the flourishing of this field did not solve all problems.
In fact, some of the most fundamental questions remained open.
Moreover, since Scheepers' pioneering work, several new notions of covers
where introduced into the framework, and some new connections with other
fields of mathematics were discovered, which helped in solving some of
the problems but introduced many new ones.

In the sequel we try to introduce a substantial portion of those problems
which lie at the core of the field. All problems presented here are interesting
enough to justify publication when solved.
However, we do not promise that all solutions will be difficult --
it could well be that we have overlooked a simple solution (there
are too many problems for us to be able to give each of them the
time it deserves). Please inform us of any solution\footnote{See
the new section \textsc{Notes added in proof} at the end of the
paper for new results obtained after the writing of the paper.}
you find or
any problem which you find important and which was not included here.
It is our hope that we will be able to publish a complementary
survey of these in the future.

Much of the material presented here is borrowed (without further notice)
from the $\mathcal{SPM}$ \textbf{Bulletin}, a semi-monthly electronic bulletin
dedicated to the field \cite{spm}.
Announcements of solutions and other problems sent to the author
will be published in this bulletin:
We urge the reader to subscribe to the $\mathcal{SPM}$ \textbf{Bulletin}\footnote{E-mail
us to get subscribed, free of charge.}
in order to remain up-to-date in this quickly evolving
field of mathematics.

We thank the organizing committee of the
\emph{Lecce Workshop on Coverings, Selections and Games in Topology}
(June 2002) for inviting this survey paper.

\subsection{A note to the reader}
The definitions always appear before the \emph{first problem}
requiring them, and are not repeated later.

\section{Basic notation and conventions}

\subsection{Selection principles}
The following notation is due to Scheepers \cite{coc1}.
This notation can be used to denote many properties which
were considered in the classical literature under various
names (see Figure~\ref{3dim} below), and
using it makes the analysis of the relationships between
these properties very convenient.

Let $\fU$ and $\fV$ be collections of covers of a space $X$.
Following are selection hypotheses which $X$ might satisfy or not
satisfy.\footnote{Often these hypotheses are identified with the class of
all spaces satisfying them.}
\begin{itemize}
\item[$\sone(\fU,\fV)$:]
For each sequence $\seq{\cU_n}$ of members of $\fU$,
there is a sequence
$\seq{U_n}$ such that for each $n$ $U_n\in\cU_n$, and $\seq{U_n}\in\fV$.
\item[$\sfin(\fU,\fV)$:]
For each sequence $\seq{\cU_n}$
of members of $\fU$, there is a sequence
$\seq{\cF_n}$ such that each $\cF_n$ is a finite
(possibly empty) subset of $\cU_n$, and
$\Union_{n\in\N}\cF_n\in\fV$.
\item[$\ufin(\fU,\fV)$:]
For each sequence
$\seq{\cU_n}$ of members of $\fU$
\emph{which do not contain a finite subcover},
there exists a sequence $\seq{\cF_n}$
such that for each $n$ $\cF_n$ is a finite (possibly empty) subset of
$\cU_n$, and
$\seq{\cup\cF_n}\in\fV$.
\end{itemize}
Clearly, $\sone(\fU,\fV)$ implies $\sfin(\fU,\fV)$, and for the
types of covers that we consider here, $\sfin(\fU,\fV)$
implies $\ufin(\fU,\fV)$.

\subsection{Stronger subcovers}
The following prototype of many classical properties is called
``\emph{$\fU$ choose $\fV$}'' in \cite{tautau}.
\bi
\i[$\binom{\fU}{\fV}$:]
For each $\cU\in\fU$ there exists $\cV\sbst\cU$ such that $\cV\in\fV$.
\ei
Then $\sfin(\fU,\fV)$ implies $\binom{\fU}{\fV}$.

\subsection{The spaces considered}

Many of the quoted results apply in the case that the spaces $X$ in question are
Tychonoff, perfectly normal, or Lindel\"of in all powers.
However, unless otherwise indicated, we consider spaces $X$ which are (homeomorphic
to) sets of reals. This is the case, e.g., for any separable zero-dimensional
metrizable space.

This significantly narrows our scope, but since we
are interested in finding \emph{good problems} rather than proving
general results, this may be viewed as a tool
to filter out problems arising from topologically-pathological examples.
However, most of the problems make sense in the general case and solutions
in the general setting are usually also of interest.

\part{Classical types of covers}

\section{Thick covers}
In this paper, by \emph{cover} we mean a nontrivial one, that is,
$\cU$ is a cover of $X$ if $X=\cup\cU$ and $X\nin\cU$.
$\cU$ is:
\be
\i A \emph{large cover} of $X$ if each $x\in X$ is contained in infinitely many members of $\cU$,
\i An \emph{$\omega$-cover} of $X$ if each finite subset of $X$ is contained in some member of $\cU$; and
\i A \emph{$\gamma$-cover} of $X$ if $\cU$ is infinite, and
each $x\in X$ belongs to all but finitely many
members of $\cU$.
\ee
The large covers and the $\w$-covers are quite old.
The term ``$\gamma$-covers'' was coined in a relatively new paper \cite{GN}
but this type of covers appears at least as early as in \cite{HURE27}.

Let $\O$, $\Lambda$, $\Omega$, and $\Gamma$
denote the collections of open covers, open large covers, $\omega$-covers,
and $\gamma$-covers of $X$, respectively.
If we assume that $X$ is a set of reals (or a separable, zero-dimensional
metric space), then we may assume that all covers in these collections
are countable \cite{GN, splittability}.
Similarly, let $\B,\BL,\BO,\BG$
be the corresponding \emph{countable Borel} covers of $X$.
Often the properties obtained by applying $\sone$, $\sfin$, or $\ufin$ to a
pair of the above families of covers are called \emph{classical selection principles}.

\section{Classification}\label{classiclassif}

The following discussion is based on \cite{coc2, CBC}.
Recall that for the types of covers which we consider,
$$\sone(\fU,\fV)\impl\sfin(\fU,\fV)\impl\ufin(\fU,\fV)\mbox{ and }\binom{\fU}{\fV}$$
and $\binom{\Lambda}{\Omega}$ does not hold for a nontrivial space $X$
\cite{coc2, strongdiags}.
This rules out several of the introduced properties as trivial.
Each of our properties is monotone decreasing in the first coordinate
and increasing in the second.
In the case of $\ufin$ note that
for any class of covers $\fV$,
$\ufin(\O,\fV)$ is equivalent to $\ufin(\Gamma,\fV)$
because given an open cover $\seq{U_n}$ we may replace
it by $\seq{\Union_{i<n} U_i}$, which is a $\gamma$-cover
(unless it contains a finite subcover).

In the three-dimensional diagram of Figure~\ref{3dim},
the double lines indicate that the two properties are equivalent.
The proof of these equivalences can be found in \cite{coc1, coc2}.

\begin{figure}[!hb]
\unitlength=.8mm
\begin{picture}(141.00,145.00)(0,0)
\put(104.00,20.00){\makebox(0,0)[cc]{$\sone(\O,\O)$}}
\put(104.00,46.00){\makebox(0,0)[cc]{$\sone(\Lambda,\O)$}}
\put(73.00,46.00){\makebox(0,0)[cc]{$\sone(\Lambda,\Lambda)$}}
\put(104.00,72.00){\makebox(0,0)[cc]{$\sone(\Omega,\O)$}}
\put(73.00,72.00){\makebox(0,0)[cc]{$\sone(\Omega,\Lambda)$}}
\put(43.00,72.00){\makebox(0,0)[cc]{$\sone(\Omega,\Omega)$}}
\put(13.00,72.00){\makebox(0,0)[cc]{$\sone(\Omega,\Gamma)$}}
\put(104.00,99.00){\makebox(0,0)[cc]{$\sone(\Gamma,\O)$}}
\put(73.00,99.00){\makebox(0,0)[cc]{$\sone(\Gamma,\Lambda)$}}
\put(43.00,99.00){\makebox(0,0)[cc]{$\sone(\Gamma,\Omega)$}}
\put(13.00,99.00){\makebox(0,0)[cc]{$\sone(\Gamma,\Gamma)$}}

\put(123.00,33.00){\makebox(0,0)[cc]{$\sfin(\O,\O)$}}
\put(123.00,59.00){\makebox(0,0)[cc]{$\sfin(\Lambda,\O)$}}
\put(89.00,59.00){\makebox(0,0)[cc]{$\sfin(\Lambda,\Lambda)$}}
\put(123.00,86.00){\makebox(0,0)[cc]{$\sfin(\Omega,\O)$}}
\put(89.00,86.00){\makebox(0,0)[cc]{$\sfin(\Omega,\Lambda)$}}
\put(58.00,86.00){\makebox(0,0)[cc]{$\sfin(\Omega,\Omega)$}}
\put(27.00,86.00){\makebox(0,0)[cc]{$\sfin(\Omega,\Gamma)$}}
\put(27.00,115.00){\makebox(0,0)[cc]{$\sfin(\Gamma,\Gamma)$}}
\put(58.00,115.00){\makebox(0,0)[cc]{$\sfin(\Gamma,\Omega)$}}
\put(89.00,115.00){\makebox(0,0)[cc]{$\sfin(\Gamma,\Lambda)$}}
\put(123.00,115.00){\makebox(0,0)[cc]{$\sfin(\Gamma,\O)$}}

\put(141.00,130.00){\makebox(0,0)[cc]{$\ufin(\Gamma,\O)$}}
\put(43.00,130.00){\makebox(0,0)[cc]{$\ufin(\Gamma,\Gamma)$}}
\put(73.00,130.00){\makebox(0,0)[cc]{$\ufin(\Gamma,\Omega)$}}
\put(104.00,130.00){\makebox(0,0)[cc]{$\ufin(\Gamma,\Lambda)$}}

\put(105.00,13.00){\makebox(0,0)[cc]{$C^{\prime\prime}$
Rothberger}}
\put(11.00,64.00){\makebox(0,0)[cc]{$\gamma$-set Gerlits-Nagy}}
\put(131.00,137.00){\makebox(0,0)[cc]{Menger}}
\put(43.00,137.00){\makebox(0,0)[cc]{Hurewicz}}

\put(054.00,130.00){\vector(1,0){9.00}}
\put(084.00,130.00){\vector(1,0){9.00}}
\put(114.00,130.50){\line(1,0){16.00}}
\put(114.00,129.50){\line(1,0){16.00}}

\put(061.00,118.00){\vector(1,1){10.00}}
\put(092.00,118.00){\line(1,1){10.00}}
\put(091.00,118.00){\line(1,1){10.00}}
\put(030.00,118.00){\vector(1,1){10.00}}
\put(122.00,118.00){\line(1,1){10.00}}
\put(123.00,118.00){\line(1,1){10.00}}

\put(036.00,115.00){\vector(1,0){12.00}}
\put(067.00,115.00){\vector(1,0){12.00}}
\put(099.00,115.50){\line(1,0){14.00}}
\put(099.00,114.50){\line(1,0){14.00}}

\put(014.50,102.00){\line(1,1){10.00}}
\put(015.50,102.00){\line(1,1){10.00}}
\put(043.00,102.00){\vector(1,1){10.00}}
\put(076.00,102.00){\vector(1,1){10.00}}
\put(107.00,102.00){\vector(1,1){10.00}}

\put(090.00,101.00){\line(0,1){10.00}}
\put(089.00,101.00){\line(0,1){10.00}}
\put(058.00,101.00){\vector(0,1){10.00}}
\put(027.00,101.00){\vector(0,1){10.00}}

\put(021.00,099.00){\vector(1,0){14.00}}
\put(051.00,099.00){\vector(1,0){14.00}}
\put(081.00,099.50){\line(1,0){14.00}}
\put(081.00,098.50){\line(1,0){14.00}}

\put(027.00,090.00){\line(0,1){6.00}}
\put(058.00,090.00){\line(0,1){6.00}}
\put(090.00,090.00){\line(0,1){6.00}}
\put(089.00,090.00){\line(0,1){6.00}}
\put(123.00,090.00){\line(0,1){22.00}}
\put(122.00,090.00){\line(0,1){22.00}}

\put(037.00,086.00){\line(1,0){4.00}}
\put(045.00,086.00){\vector(1,0){3.00}}
\put(068.00,086.00){\line(1,0){3.00}}
\put(074.00,086.00){\vector(1,0){5.00}}
\put(099.00,086.50){\line(1,0){3.00}}
\put(099.00,085.50){\line(1,0){3.00}}
\put(106.00,086.50){\line(1,0){7.00}}
\put(106.00,085.50){\line(1,0){7.00}}

\put(016.00,076.00){\line(1,1){7.00}}
\put(017.00,076.00){\line(1,1){7.00}}
\put(047.00,076.00){\vector(1,1){7.00}}
\put(077.00,076.00){\vector(1,1){7.00}}
\put(107.00,076.00){\vector(1,1){7.00}}

\put(013.00,076.00){\vector(0,1){19.00}}
\put(043.00,076.00){\vector(0,1){19.00}}
\put(072.00,076.00){\vector(0,1){20.00}}
\put(090.00,075.00){\line(0,1){7.00}}
\put(089.00,075.00){\line(0,1){7.00}}
\put(104.00,076.00){\vector(0,1){20.00}}

\put(021.00,072.00){\vector(1,0){14.00}}
\put(051.00,072.00){\vector(1,0){14.00}}
\put(081.00,072.50){\line(1,0){14.00}}
\put(081.00,071.50){\line(1,0){14.00}}

\put(089.00,062.00){\line(0,1){7.00}}
\put(090.00,062.00){\line(0,1){7.00}}
\put(122.00,064.00){\line(0,1){19.00}}
\put(123.00,064.00){\line(0,1){19.00}}

\put(099.00,059.50){\line(1,0){3.00}}
\put(099.00,058.50){\line(1,0){3.00}}
\put(106.00,059.50){\line(1,0){6.00}}
\put(106.00,058.50){\line(1,0){6.00}}

\put(078.00,050.00){\vector(1,1){7.00}}
\put(108.00,050.00){\vector(1,1){7.00}}

\put(073.00,050.00){\line(0,1){19.00}}
\put(072.00,050.00){\line(0,1){19.00}}
\put(104.00,049.00){\line(0,1){19.00}}
\put(103.00,049.00){\line(0,1){19.00}}

\put(082.00,046.00){\line(1,0){13.00}}
\put(082.00,045.00){\line(1,0){13.00}}

\put(123.00,035.00){\line(0,1){19.00}}
\put(122.00,035.00){\line(0,1){19.00}}

\put(104.00,023.00){\line(0,1){19.00}}
\put(103.00,023.00){\line(0,1){19.00}}

\put(108.00,023.00){\vector(1,1){8.00}}

\end{picture}
\caption{\label{3dim}}
\end{figure}
The analogous equivalences for the Borel case also hold, but
in the Borel case more equivalences hold \cite{CBC}:
For each $\fV \in\{\B,\BO,\BG\}$,
$$\sone(\BG,\fV) \Iff \sfin(\BG,\fV) \Iff \ufin(\BG,\fV).$$
After removing duplications we obtain Figure \ref{survive}.

\begin{figure}[!htp]
\newcommand{\sr}[2]{#1}
{\tiny
$$\xymatrix@C=-2pt@R=15pt{
&
&
& \sr{\ufin(\Gamma,\Gamma)}{\b}\ar[rr]\ar@{.>}[dr]^?
&
& \sr{\ufin(\Gamma,\Omega)}{\d}\ar[rrrrr]\ar@/_/@{.>}[dl]_?
&
&
&
&
&
&
& \sr{\ufin(\Gamma,\O)}{\d}
\\
&
&
&
& \sr{\sfin(\Gamma,\Omega)}{\d}\ar[ur]
\\
& \sr{\sone(\Gamma,\Gamma)}{\b}\ar[rr]\ar[uurr]
&
& \sr{\sone(\Gamma,\Omega)}{\d}\ar[rrr]\ar[ur]
&
&
& \sr{\sone(\Gamma,\O)}{\d}\ar[uurrrrrr]
\\
  \sr{\sone(\BG,\BG)}{\b}\ar[ur]\ar[rr]
&
& \sr{\sone(\BG,\BO)}{\d}\ar[ur]\ar[rrr]
&
&
& \sr{\sone(\BG,\B)}{\d}\ar[ur]
\\
&
&
&
& \sr{\sfin(\Omega,\Omega)}{\d}\ar'[u]'[uu][uuu]
\\
\\
&
& \sr{\sfin(\BO,\BO)}{\d}\ar[uuu]\ar[uurr]
\\
& \sr{\sone(\Omega,\Gamma)}{\p}\ar'[r][rr]\ar'[uuuu][uuuuu]
&
& \sr{\sone(\Omega,\Omega)}{\cov(\M)}\ar'[uuuu][uuuuu]\ar'[rr][rrr]\ar[uuur]
&
&
& \sr{\sone(\O,\O)}{\cov(\M)}\ar[uuuuu]
\\
  \sr{\sone(\BO,\BG)}{\p}\ar[uuuuu]\ar[rr]\ar[ur]
&
& \sr{\sone(\BO,\BO)}{\cov(\M)}\ar[uu]\ar[ur]\ar[rrr]
&
&
& \sr{\sone(\B,\B)}{\cov(\M)}\ar[uuuuu]\ar[ur]
}$$
}
\caption{}\label{survive}
\end{figure}

All implications which do not appear in Figure \ref{survive} where
refuted by counter-examples (which are in fact sets of real numbers)
in \cite{coc1, coc2, CBC}.
The only unsettled implications in this diagram are marked with
dotted arrows.

\begin{prob}[{\Cite{Problems 1 and 2}{coc2}}]\label{finishclass}
~\be
\i Is $\ufin(\Gamma,\Omega)$ equivalent to $\sfin(\Gamma,\Omega)$?
\i And if not, does $\ufin(\Gamma,\Gamma)$ imply $\sfin(\Gamma,\Omega)$?
\ee
\end{prob}
Bartoszy\'nski (personal communication) suspects
that an implication should be easy to prove if it is true,
and otherwise it may be quite difficult to find a counter-example
(existing methods do not tell these properties apart).
However, the Hurewicz property $\ufin(\Gamma,\Gamma)$ has some surprising disguises
which a priori do not look equivalent to it \cite{coc7, hureslaloms},
so no definite conjecture can be made about this problem.

\section{Classification in ZFC}
Most of the examples used to prove non-implications in Figure~\ref{survive}
are ones using (fragments of) \CH{}.
However, some non-implications can be proved without any extra hypotheses.
For example, every $\sigma$-compact space satisfies
$\ufin(\Gamma,\Gamma)$ and $\sfin(\Omega,\Omega)$ (and all properties implied by these),
but the Cantor set does not satisfy $\sone(\Gamma,\O)$ (and all properties implying it) \cite{coc2}.

It is not known if additional non-implications are provable without the help of additional
axioms. We mention one problem which drew more attention then the others.

\begin{prob}[\Cite{Problem 3}{coc2}, \Cite{Problem 1}{BuHa}, \Cite{Problem 1}{LBopenprobs}]\label{MnotH}
Does there exist (in ZFC) a set of reals $X$ which has the Menger property
$\ufin(\O,\O)$ but not the Hurewicz property $\ufin(\O,\Gamma)$?
\end{prob}

Not much is known about the situation when arbitrary topological spaces
are considered rather than sets of reals.
\begin{proj}[{\cite[Problem 3]{coc2}}]\label{generaltop}
Find, without extra hypotheses, (general) topological spaces that demonstrate non-implications among the
classical properties. Do the same for \emph{Lindel\"of} topological spaces.
\end{proj}

\section{Uncountable elements in ZFC}\label{ZFCelements}

We already mentioned the fact that the Cantor set satisfies all
properties of type $\sfin$ or $\ufin$ in the case of open covers,
but none of the remaining ones. It turns out that some $\sone$-type
properties can be shown to be satisfied by uncountable elements without
any special hypotheses.

This is intimately related to the following notions.
The \emph{Baire space} $\NN$ is equipped with the product topology
and (quasi)ordered by eventual dominance: $f\le^* g$ if
$f(n)\le^* g(n)$ for all but finitely many $n$.
A subset of $\NN$ is \emph{dominating} if it is cofinal in $\NN$ with respect to $\le^*$.
If a subset of $\NN$ is unbounded with respect to $\le^*$ then we simply say that it is
\emph{unbounded}.
Let $\b$ (respectively, $\d$) denote the minimal cardinality of an unbounded
(respectively, dominating) subset of $\NN$.

The \emph{critical cardinality} of a nontrivial family $\cJ$ of sets of reals
is $\non(\cJ)=\min\{|X| : X\sbst\R\mbox{ and }X\nin\cJ\}$.
Then $\b$ is the critical cardinality of $\sone(\BG,\BG)$, $\sone(\Gamma,\Gamma)$, and
$\ufin(\Gamma,\Gamma)$, and $\d$ is the critical cardinality of the classes in
Figure~\ref{survive} which contain $\sfin(\BO,\BO)$ \cite{coc2, CBC}.

\subsection{The open case}
In \cite{coc2, wqn} it was shown (in ZFC) that there exists a set of reals
of size $\aleph_1$ which satisfies $\sone(\Gamma,\Gamma)$ as well as $\sfin(\Omega,\Omega)$.
In \cite{alpha_i} this is improved to show that there always exists a set of
size $\t$ which satisfies these properties. In both cases the proof uses a dichotomy
argument (two different examples are given in two possible extensions of ZFC).

In \cite{ideals} the following absolute ZFC example is studied.
Let $\compactN$ be the one point compactification of $\N$.
(A subset $A\sbst\compactN$ is open if: $A\sbst\N$, or $\infty\in A$ and $A$ is cofinite.)
Let $\Increasing \sbst \NcompactN$ consist of the functions $f$ such
that
\be
\item For all $n$, $f(n) \leq f(n+1)$; and
\item For all $n$, if $f(n)<\infty$, then $f(n)<f(n+1)$.
\ee
($\Increasing$ is homeomorphic to the Cantor set of reals.)
For each increasing finite sequence $s$ of natural numbers,
let $q_{s} \in \Increasing$ be defined as
$$q_{s}(k)=
\begin{cases}
s(k) & \text{if } k <|s|\\
\infty & \text{otherwise}
\end{cases}$$
for each $k\in\N$.
Note that the set
$$Q=\{q_s : s\text{ an increasing finite sequence in }\N \}$$
is dense in $\Increasing$.

Let $B=\{f_\alpha : \alpha < \b\}\sbst\NN$ be a $\le^*$-unbounded
set of strictly increasing elements of $\NN$ which forms a \emph{$\b$-scale}
(that is, for each $\alpha<\beta$, $f_\alpha \leq^* f_\beta$),
and set $H = B\cup Q$.
In \cite{ideals} it is proved that
all finite powers of $H$ satisfy $\ufin(\O,\Gamma)$.

\begin{prob}[\Cite{Problem 17}{ideals}, \Cite{Problem 2}{o-bdd}]\label{isS1GG}
Does $H$ satisfy $\sone(\Gamma,\Gamma)$?
\end{prob}
By the methods of \cite{wqn}, it would be enough to prove that
\begin{quote}
For each sequence $\seq{\cU_n}$ of open $\gamma$-covers
of $X$, there exists a sequence $\seq{U_n}$ such that
for each $n$ $U_n\in\cU_n$, and a subset $Y\sbst X$,
such that $|Y|<\b$ and $\seq{U_n}$ is a $\gamma$-cover
of $X\sm Y$.
\end{quote}
to obtain a positive answer.

\subsection{The Borel case}
Borel's Conjecture, which was proved to be consistent by Laver,
implies that each set of reals satisfying $\sone(\O,\O)$ (and the classes below it) is
countable.
From our point of view this means that there do not exist ZFC examples
of sets satisfying $\sone(\O,\O)$.
A set of reals $X$ is a \emph{$\sigma$-set} if each $G_\delta$ set in $X$
is also an $F_\sigma$ set in $X$.
In \cite{CBC} it is proved that every element of $\sone(\BG,\BG)$
is a $\sigma$-set.
According to a result of Miller \cite{Mil79Len},
it is consistent that every $\sigma$-set of
real numbers is countable. Thus, there do not exist uncountable
ZFC examples satisfying $\sone(\BG,\BG)$.
The situation for the other classes, though addressed by top
experts, remains open. In particular, we have the following.

\begin{prob}[\cite{pawlikowskireclaw}, \Cite{Problem 45}{CBC}, \cite{ideals}]
Does there exist (in ZFC) an uncountable set of reals satisfying $\sone(\BG,\B)$?
\end{prob}
By \cite{CBC}, this is the same as asking whether it is consistent
that each uncountable set of reals can be mapped onto a dominating
subset of $\NN$ by a Borel function.
This is one of the major open problems in the field.

\section{Special elements under weak hypotheses}

Most of the counter examples used to distinguish between the properties
in the Borel case are constructed with the aid of \CH{}.
The question whether such examples exist under weaker hypotheses
(like \MA{}) is often raised (e.g., \cite{JORG, MilNonGamma}).
We mention some known results
by indicating (by full bullets) all places in the diagram of the
Borel case (the front plane in Figure~\ref{survive})
which the example satisfies.
All hypotheses we mention are weaker than \MA{}.

Let us recall the basic terminology.
$\M$ and $\cN$ denote the collections of meager (=first category)
and null (=Lebesgue measure zero) sets of reals, respectively.
For a family $\cI$ of sets of reals define:
\begin{eqnarray*}
\add(\cI) & = & \min\{|\cF| : \cF\sbst\cI\mbox{ and } \cup\cF\nin\cI\}\\
\cov(\cI) & = & \min\{|\cF| : \cF\sbst\cI\mbox{ and } \cup\cF=\R\}\\
\cof(\cI) & = & \min\{|\cF| : \cF\sbst\cI\mbox{ and } (\forall I\in\cI)(\exists J\in\cF)\ I\sbst J\}
\end{eqnarray*}
A set of reals $X$ is a \emph{$\kappa$-Luzin} set if
$|X|\ge\kappa$ and for each $M\in\M$, $|X\cap M|<\kappa$.
Dually, $X$ is a \emph{$\kappa$-Sierpi\'nski} set
if $|X|\ge\kappa$ and for each $N\in\cN$, $|X\cap N|<\kappa$.

\begin{figure}[!ht]
\begin{center}
$\xymatrix@C=10pt@R=12pt{
\circ\ar[r]        & \circ\ar[r]       & \bullet        \\
                    & \circ\ar[u]                        \\
\circ\ar[uu]\ar[r] & \circ\ar[u]\ar[r] & \bullet\ar[uu]\\
                    & (a)
}$
\quad
$\xymatrix@C=10pt@R=12pt{
\circ \ar[r]        & \bullet \ar[r]       & \bullet        \\
                    & \bullet \ar[u]                        \\
\circ \ar[uu]\ar[r] & \bullet \ar[u]\ar[r] & \bullet \ar[uu]\\
                    & (b)
}$
\quad
$\xymatrix@C=10pt@R=12pt{
\bullet \ar[r]        & \bullet \ar[r]       & \bullet        \\
                      & \circ   \ar[u]                        \\
\circ   \ar[uu]\ar[r] & \circ   \ar[u]\ar[r] & \circ   \ar[uu]\\
                      & (c)
}$
\quad
$\xymatrix@C=10pt@R=12pt{
\bullet \ar[r]        & \bullet \ar[r]       & \bullet        \\
                      & \bullet \ar[u]                        \\
\circ   \ar[uu]\ar[r] & \circ   \ar[u]\ar[r] & \circ   \ar[uu]\\
                      & (d)
}$
\end{center}
\caption{Some known configurations under \MA{}}\label{bullet}
\end{figure}
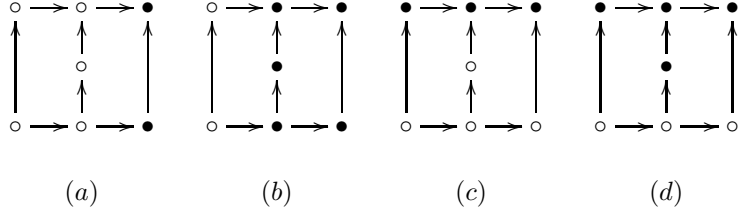

If $\cov(\M)=\cof(\M)$ then there exists a $\cov(\M)$-Luzin set satisfying
the properties indicated in Figure~\ref{bullet}(a) \cite{CBC}
(in \cite{luzinundetermined} this is proved under \MA{}).
Under the slightly stronger assumption $\cov(\M)=\c$,
there exists a $\cov(\M)$-Luzin set as in Figure~\ref{bullet}(b) \cite{huremen2}.\footnote{%
In fact, we can require that this Luzin set does not satisfy $\ufin(\BG,\BT)$
\cite{tautau} -- see Section \ref{tau} for the definition of
$\BT$.} Dually, assuming $\cov(\cN)=\cof(\cN)=\b$ there exists a
$\b$-Sierpi\'nski set
as in Figure~\ref{bullet}(c) \cite{CBC}, and another one as in Figure~\ref{bullet}(d) \cite{o-bdd}.\footnote{%
The last Sierpi\'nski set actually satisfies $\sfin(\BO^\gpbl,\BO)$ -- see Section \ref{groupbl}
for the definition of $\BO^\gpbl$.
}

\forget 
%
\forgotten

\begin{proj}
Find constructions, under \MA{} or weaker hypotheses, for any of
the consistent configurations not covered in Figure~\ref{bullet}.
\end{proj}

\section{Preservation of properties}

\subsection{Hereditariness}
A property is (provably) \emph{hereditary} if for each space $X$ satisfying the property,
all subsets of $X$ satisfy that property.
Most of the properties considered in this paper may be considered intuitively
as notions of smallness, thus it is somewhat surprising that none of the properties involving
open covers is hereditary \cite{ideals}.
However, the property $\sone(\B,\B)$ as well as all properties of the form $\Pi(\BG,\B)$
are hereditary \cite{ideals} (but $\sone(\BO,\BG)$ is not \cite{MilNonGamma}).

\begin{prob}[\Cite{Problem 4}{ideals}, \Cite{Question 6}{MilNonGamma}]\label{borhered}
Is any of the properties $\sone(\BO,\BO)$ or $\sfin(\BO,\BO)$ hereditary?
\end{prob}

This problem is related to Problem \ref{borpows} below.

\subsection{Finite powers}
$\sone(\Omega,\Gamma)$, $\sone(\Omega,\Omega)$, and $\sfin(\Omega,\Omega)$ are the
only properties in the open case which are preserved under taking finite powers
\cite{coc2}.
The only candidates in the Borel case to be preserved under taking finite powers
are the following.

\begin{prob}[\Cite{Problem 50}{CBC}]\label{borpows1}
Is any of the classes $\sone(\BO,\BG)$, $\sone(\BO,\BO)$, or $\sfin(\BO,\BO)$
closed under taking finite powers?
\end{prob}

In \cite{CBC} it is shown that if all finite powers of $X$ satisfy $\sone(\B,\B)$
(respectively, $\sone(\BG,\B)$), then $X$ satisfies $\sone(\BO,\BO)$
(respectively, $\sfin(\BO,\BO)$). Consequently, the last two cases of Problem \ref{borpows1}
translate to the following.

\begin{prob}[\Cite{Problems 19 and 21}{CBC}]\label{borpows}
~\be
\i Is it true that if $X$ satisfies $\sone(\BO,\BO)$, then all finite powers of $X$ satisfy $\sone(\B,\B)$?
\i Is it true that if $X$ satisfies $\sfin(\BO,\BO)$, then all finite powers of $X$ satisfy $\sone(\BG,\B)$?
\ee
\end{prob}
The analogous assertion in the open case is true \cite{sakai, coc2}.
Observe that a positive answer to this problem implies a positive answer to Problem \ref{borhered} above.

It is worthwhile to mention that by a sequence of results in \cite{lengthdiags, CBC, huremen2, splittability},
none of the properties in Figure~\ref{survive} is preserved under taking finite products.

\subsection{Unions}
The question of which of the properties in Figure~\ref{survive} is provably
preserved under taking finite or countable unions (that is, finitely or countably additive)
is completely settled in \cite{coc2, wqn, lengthdiags, huremen2}.
Also, among the classes which are not provably additive, it is known that some are
\emph{consistently} additive \cite{huremen2}.
Only the following problems remain open in this category.

\begin{prob}[\cite{AddQuad}]
Is any of the classes
$\sfin(\Omega,\Omega)$, $\sone(\Gamma,\Omega)$, and $\sfin(\Gamma,\Omega)$
consistently closed under taking finite unions?
\end{prob}

\begin{prob}[\cite{AddQuad}]
Is $\sfin(\BO,\BO)$ consistently closed under taking finite unions?
\end{prob}

Another sort of problems which remain open is that of
determining the \emph{exact} additivity numbers of the (provably)
additive properties.
The general problem is to determine the additivity numbers of the
properties in question in terms of well known cardinal characteristics
of the continuum (like $\b$, $\d$, etc.). See \cite{LecceSurvey} for
a list of properties for which the problem is still open.
Three of the more interesting ones among these are the following.

\begin{prob}[\cite{AddQuad}]\label{addM}
Is $\add(\ufin(\Gamma,\O))=\b$?
\end{prob}
It is only known that $\b\le\add(\ufin(\Gamma,\O))\le\cf(\d)$, and
that the additivity of the corresponding combinatorial notion of
smallness is equal to $\b$ \cite{AddQuad}.

\begin{prob}\label{addGGb}
Is $\add(\sone(\Gamma,\Gamma))=\b$?
\end{prob}
It is known that $\h\le\add(\sone(\Gamma,\Gamma))\le\b$ \cite{wqn}.
This problem is related to Problem \ref{s1add} below.

In \cite{covM2} it is proved that $\add(\cN)\le\add(\sone(\O,\O))$.
\begin{prob}[\Cite{Problem 4}{covM2}]
Is it consistent that $\add(\cN)<\add(\sone(\O,\O))$?
\end{prob}

\part{Modern types of covers}

In this part we divide the problems according to the involved
type of covers rather than according to the type of problem.

\section{$\tau$-covers}\label{tau}

$\cU$ is a \emph{$\tau$-cover} of $X$ \cite{tau}
if it is a large cover of $X$,\footnote{Recall
that by ``cover of $X$'' we mean one not containing $X$ as an element.}
and for each $x,y\in X$, either $\{U\in\cU : x\in U, y\nin U\}$ is finite, or
else $\{U\in\cU : y\in U, x\nin U\}$ is finite.
If all powers of $X$ are Lindel\"of (e.g, if $X$ is a set of reals) then
each $\tau$-cover of $X$ contains a countable $\tau$-cover of $X$
\cite{splittability}.
Let $\Tau$ denote the collection of open $\tau$-covers of $X$. Then
$$\Gamma  \sbst \Tau \sbst \Omega \sbst \O.$$

The following problem arises in almost every study of $\tau$-covers
\cite{tau, tautau, splittability, strongdiags, Hdim, MilNonGamma}.
By \cite{GN}, $\sone(\Omega,\Gamma)\Iff\binom{\Omega}{\Gamma}$.
As $\Gamma\sbst\Tau$, this property implies $\binom{\Omega}{\Tau}$.
Thus far, all examples of sets not satisfying $\binom{\Omega}{\Gamma}$
turned out not to satisfy $\binom{\Omega}{\Tau}$.
\begin{prob}[\Cite{\S 4}{spm1}]\label{omtau}
Is $\binom{\Omega}{\Gamma}$ equivalent to $\binom{\Omega}{\Tau}$?
\end{prob}
A positive answer would imply that the properties $\sone(\Omega,\Gamma)$, $\sone(\Omega,\Tau)$,
and $\sfin(\Omega,\Tau)$ are all equivalent, and therefore simplify the study of $\tau$-covers
considerably. It would also imply a positive solution to Problems \ref{schtauprob}, \ref{tauprod}(1),
\ref{splittabilityclassif}, and other problems.\footnote{This looks too good to be true,
but a negative answer should also imply (through a bit finer analysis) a solution to several
open problems.}
The best known result in this direction is that $\binom{\Omega}{\Tau}$
implies $\sfin(\Gamma,\Tau)$ \cite{tautau}.
A modest form of Problem \ref{omtau} is the following.
If $\binom{\Omega}{\Tau}$ implies $\sfin(\Tau,\Omega)$, then $\binom{\Omega}{\Tau}\Iff\sfin(\Omega,\Tau)$.
\begin{prob}[\Cite{Problem 2.9}{tautau}]
Is $\binom{\Omega}{\Tau}$ equivalent to $\sfin(\Omega,\Tau)$?
\end{prob}

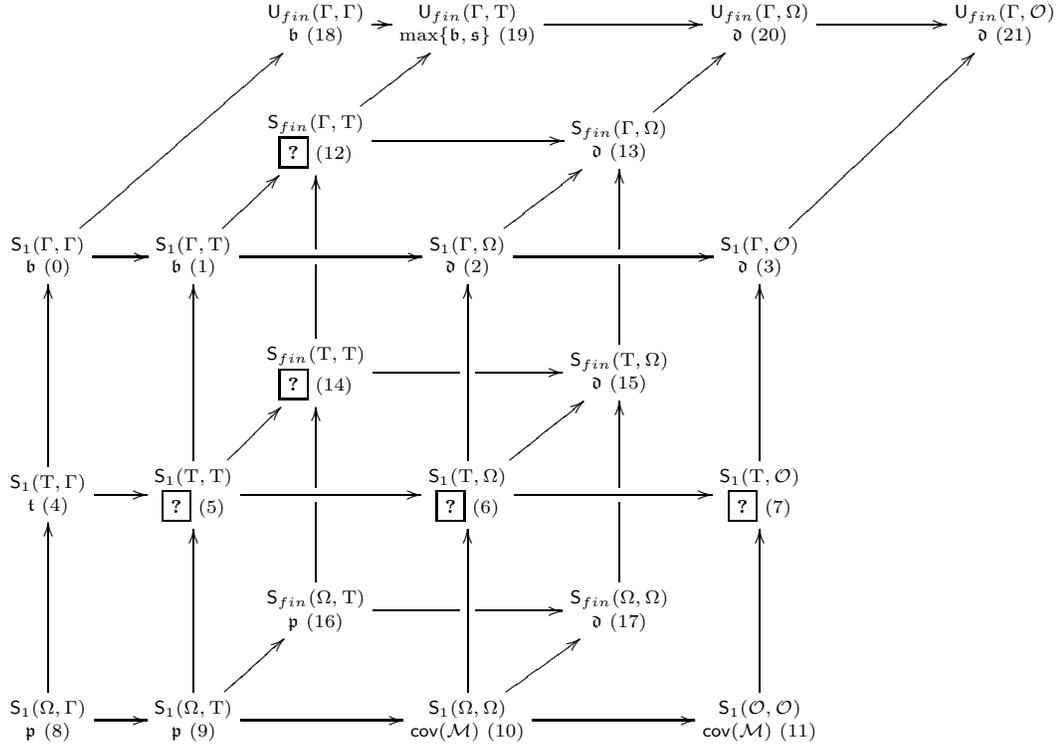
\begin{figure}[!ht]
\newcommand{\sr}[2]{{\txt{$#1$\\$#2$}}}
{\scriptsize
\begin{center}
$\xymatrix@C=7pt@R=20pt{
&
&
& \sr{\ufin(\Gamma,\Gamma)}{\b~~ (18)}\ar[r]
& \sr{\ufin(\Gamma,\Tau)}{\max\{\b,\s\}~~ (19)}\ar[rr]
&
& \sr{\ufin(\Gamma,\Omega)}{\d~~ (20)}\ar[rrrr]
&
&
&
& \sr{\ufin(\Gamma,\O)}{\d~~ (21)}
\\
&
&
& \sr{\sfin(\Gamma,\Tau)}{\fbox{\textbf{?}}~~ (12)}\ar[rr]\ar[ur]
&
& \sr{\sfin(\Gamma,\Omega)}{\d~~ (13)}\ar[ur]
\\
\sr{\sone(\Gamma,\Gamma)}{\b~~ (0)}\ar[uurrr]\ar[rr]
&
& \sr{\sone(\Gamma,\Tau)}{\b~~ (1)}\ar[ur]\ar[rr]
&
& \sr{\sone(\Gamma,\Omega)}{\d~~ (2)}\ar[ur]\ar[rr]
&
& \sr{\sone(\Gamma,\O)}{\d~~ (3)}\ar[uurrrr]
\\
&
&
& \sr{\sfin(\Tau,\Tau)}{\fbox{\textbf{?}}~~ (14)}\ar'[r][rr]\ar'[u][uu]
&
& \sr{\sfin(\Tau,\Omega)}{\d~~ (15)}\ar'[u][uu]
\\
\sr{\sone(\Tau,\Gamma)}{\t~~ (4)}\ar[rr]\ar[uu]
&
& \sr{\sone(\Tau,\Tau)}{\fbox{\textbf{?}}~~ (5)}\ar[uu]\ar[ur]\ar[rr]
&
& \sr{\sone(\Tau,\Omega)}{\fbox{\textbf{?}}~~ (6)}\ar[uu]\ar[ur]\ar[rr]
&
& \sr{\sone(\Tau,\O)}{\fbox{\textbf{?}}~~ (7)}\ar[uu]
\\
&
&
& \sr{\sfin(\Omega,\Tau)}{\p~~ (16)}\ar'[u][uu]\ar'[r][rr]
&
& \sr{\sfin(\Omega,\Omega)}{\d~~ (17)}\ar'[u][uu]
\\
\sr{\sone(\Omega,\Gamma)}{\p~~ (8)}\ar[uu]\ar[rr]
&
& \sr{\sone(\Omega,\Tau)}{\p~~ (9)}\ar[uu]\ar[ur]\ar[rr]
&
& \sr{\sone(\Omega,\Omega)}{\cov(\M)~~ (10)}\ar[uu]\ar[ur]\ar[rr]
&
& \sr{\sone(\O,\O)}{\cov(\M)~~ (11)}\ar[uu]
}$
\end{center}
}

\caption{The diagram involving $\tau$-covers}\label{tausurvive}
\end{figure}

The notion of $\tau$-covers introduces seven new pairs---namely,
$(\Tau,\O)$, $(\Tau,\Omega)$, $(\Tau,\Tau)$, $(\Tau,\Gamma)$,
$(\O,\Tau)$, $(\Omega,\Tau)$, and $(\Gamma,\Tau)$---to
which any of the selection operators $\sone$, $\sfin$, and $\ufin$ can be applied.
This makes a total of $21$ new selection hypotheses.
Fortunately, some of them are easily eliminated.
The surviving properties appear in Figure \ref{tausurvive}.

Below each property in Figure \ref{tausurvive} appears a ``serial number''
(to be used in Table \ref{impltab}), and its critical cardinality.
The cardinal numbers $\p$, $\t$, and $\s$
are the well-known pseudo-intersection number, tower number, and splitting
number (see, e.g., \cite{vD} or \cite{BlassHBK} for definitions and details).

As indicated in the diagram, some of the critical cardinalities are not yet known.
\begin{proj}[\Cite{Problem 6.6}{tautau}]\label{cards}
What are the unknown critical cardinalities in Figure \ref{tausurvive}?
\end{proj}
Recall that there are only two unsettled implications in
the corresponding diagram for the \emph{classical} types of open covers
(Section \ref{classiclassif}).
As there are many more properties when $\tau$-covers are incorporated into the
framework, and since this investigation is new, there remain \emph{many} unsettled
implications in Figure \ref{tausurvive}. To be precise, there are exactly $76$ unsettled
implications in this diagram. These appear as question marks in the
\emph{Implications Table} \ref{impltab}.
Entry $(i,j)$ in the table ($i$th row, $j$th column) is to be interpreted as follows:
It is $1$ if property $i$ implies property $j$, $0$ if property $i$ does not
imply property $j$ (that is, consistently there exists a counter-example),
and $?$ if the implication is unsettled.

\newcommand{\mb}[1]{{\mbox{\textbf{#1}}}}

\begin{table}[!ht]
\begin{changemargin}{-5cm}{-5cm}
\begin{center}
{\scriptsize
\begin{tabular}{|r||cccccccccccccccccccccc|}
\hline
   & \mb{0} & \mb{1} & \mb{2} & \mb{3} & \mb{4} & \mb{5} & \mb{6} & \mb{7} & \mb{8} & \mb{9} & \mb{10} & \mb{11} & \mb{12} & \mb{13} & \mb{14} & \mb{15} & \mb{16} & \mb{17} & \mb{18} & \mb{19} & \mb{20} & \mb{21}\cr
\hline\hline
\mb{ 0} &  1 & 1 & 1 & 1 & 0 & \mb{?} & \mb{?} & \mb{?} & 0 & 0 & 0 & 0 & 1 & 1 & \mb{?} & \mb{?} & 0 & 0 & 1 & 1 & 1 & 1\cr
\mb{ 1} &  \mb{?} & 1 & 1 & 1 & 0 & \mb{?} & \mb{?} & \mb{?} & 0 & 0 & 0 & 0 & 1 & 1 & \mb{?} & \mb{?} & 0 & 0 & \mb{?} & 1 & 1 & 1 \cr
\mb{ 2} &  0 & 0 & 1 & 1 & 0 & 0 & \mb{?} & \mb{?} & 0 & 0 & 0 & 0 & 0 & 1 & 0 & \mb{?} & 0 & 0 & 0 & 0 & 1 & 1 \cr
\mb{ 3} &  0 & 0 & 0 & 1 & 0 & 0 & 0 & \mb{?} & 0 & 0 & 0 & 0 & 0 & 0 & 0 & 0 & 0 & 0 & 0 & 0 & 0 & 1 \cr
\mb{ 4} &  1 & 1 & 1 & 1 & 1 & 1 & 1 & 1 & 0 & 0 & \mb{?} & \mb{?} & 1 & 1 & 1 & 1 & 0 & \mb{?} & 1 & 1 & 1 & 1 \cr
\mb{ 5} &  \mb{?} & 1 & 1 & 1 & \mb{?} & 1 & 1 & 1 & 0 & 0 & \mb{?} & \mb{?} & 1 & 1 & 1 & 1 & 0 & \mb{?} & \mb{?} & 1 & 1 & 1 \cr
\mb{ 6} &  0 & 0 & 1 & 1 & 0 & 0 & 1 & 1 & 0 & 0 & \mb{?} & \mb{?} & 0 & 1 & 0 & 1 & 0 & \mb{?} & 0 & 0 & 1 & 1 \cr
\mb{ 7} &  0 & 0 & 0 & 1 & 0 & 0 & 0 & 1 & 0 & 0 & 0 & \mb{?} & 0 & 0 & 0 & 0 & 0 & 0 & 0 & 0 & 0 & 1 \cr
\mb{ 8} &  1 & 1 & 1 & 1 & 1 & 1 & 1 & 1 & 1 & 1 & 1 & 1 & 1 & 1 & 1 & 1 & 1 & 1 & 1 & 1 & 1 & 1 \cr
\mb{ 9} &  \mb{?} & 1 & 1 & 1 & \mb{?} & 1 & 1 & 1 & \mb{?} & 1 & 1 & 1 & 1 & 1 & 1 & 1 & 1 & 1 & \mb{?} & 1 & 1 & 1 \cr
\mb{10} &  0 & 0 & 1 & 1 & 0 & 0 & 1 & 1 & 0 & 0 & 1 & 1 & 0 & 1 & 0 & 1 & 0 & 1 & 0 & 0 & 1 & 1 \cr
\mb{11} &  0 & 0 & 0 & 1 & 0 & 0 & 0 & 1 & 0 & 0 & 0 & 1 & 0 & 0 & 0 & 0 & 0 & 0 & 0 & 0 & 0 & 1 \cr
\mb{12} &  \mb{?} & \mb{?} & \mb{?} & \mb{?} & 0 & \mb{?} & \mb{?} & \mb{?} & 0 & 0 & 0 & 0 & 1 & 1 & \mb{?} & \mb{?} & 0 & 0 & \mb{?} & 1 & 1 & 1 \cr
\mb{13} &  0 & 0 & 0 & 0 & 0 & 0 & 0 & 0 & 0 & 0 & 0 & 0 & 0 & 1 & 0 & \mb{?} & 0 & 0 & 0 & 0 & 1 & 1 \cr
\mb{14} &  \mb{?} & \mb{?} & \mb{?} & \mb{?} & \mb{?} & \mb{?} & \mb{?} & \mb{?} & 0 & 0 & \mb{?} & \mb{?} & 1 & 1 & 1 & 1 & 0 & \mb{?} & \mb{?} & 1 & 1 & 1 \cr
\mb{15} &  0 & 0 & 0 & 0 & 0 & 0 & 0 & 0 & 0 & 0 & 0 & 0 & 0 & 1 & 0 & 1 & 0 & \mb{?} & 0 & 0 & 1 & 1 \cr
\mb{16} &  \mb{?} & \mb{?} & \mb{?} & \mb{?} & \mb{?} & \mb{?} & \mb{?} & \mb{?} & \mb{?} & \mb{?} & \mb{?} & \mb{?} & 1 & 1 & 1 & 1 & 1 & 1 & \mb{?} & 1 & 1 & 1 \cr
\mb{17} &  0 & 0 & 0 & 0 & 0 & 0 & 0 & 0 & 0 & 0 & 0 & 0 & 0 & 1 & 0 & 1 & 0 & 1 & 0 & 0 & 1 & 1 \cr
\mb{18} &  0 & 0 & 0 & 0 & 0 & 0 & 0 & 0 & 0 & 0 & 0 & 0 & 0 & \mb{?} & 0 & \mb{?} & 0 & 0 & 1 & 1 & 1 & 1 \cr
\mb{19} &  0 & 0 & 0 & 0 & 0 & 0 & 0 & 0 & 0 & 0 & 0 & 0 & 0 & \mb{?} & 0 & \mb{?} & 0 & 0 & 0 & 1 & 1 & 1 \cr
\mb{20} &  0 & 0 & 0 & 0 & 0 & 0 & 0 & 0 & 0 & 0 & 0 & 0 & 0 & \mb{?} & 0 & \mb{?} & 0 & 0 & 0 & 0 & 1 & 1 \cr
\mb{21} &  0 & 0 & 0 & 0 & 0 & 0 & 0 & 0 & 0 & 0 & 0 & 0 & 0 & 0 & 0 & 0 & 0 & 0 & 0 & 0 & 0 & 1 \cr
\hline
\end{tabular}
}
\caption{Implications and nonimplications}\label{impltab}
\end{center}
\end{changemargin}
\end{table}

\begin{proj}[\Cite{Problem 6.5}{tautau}]\label{implications}
Settle any of the unsettled implications in Table \ref{impltab}.
\end{proj}
Marion Scheepers asked us which single
solution would imply as many other solutions
as possible. The answer found by a computer program is the following:
If entry $(12,5)$ is $1$ (that is, $\sfin(\Gamma,\Tau)$ implies $\sone(\Tau,\Tau)$), then
there remain only 33 (!) open problems.
The best possible negative entry is $(16,3)$: If $\sfin(\Omega,\Tau)$ does not imply
$\sone(\Gamma,\O)$, then only $47$ implications remain unsettled.

Finally, observe that any solution in Problem \ref{cards}
may imply several new nonimplications.

Scheepers chose the following problem out of all the problems discussed above
as the most interesting.
\begin{prob}\label{schtauprob}
Does $\sone(\Omega,\Tau)$ imply the Hurewicz property $\ufin(\Gamma,\Gamma)$?
\end{prob}

The reason for this choice is that
if the answer is positive, then $\sone(\Omega,\Tau)$ implies
the Gerlitz-Nagy $(*)$ property \cite{GN}, which is equivalent to
another modern selection property (see Section \ref{groupbl} below).

Not much is known about the preservation of the new properties under set theoretic
operations.
Miller \cite{MilNonGamma} proved that assuming \CH{}, there exists
a set $X$ satisfying $\sone(\BO,\BG)$ and a subset $Y$ of $X$ such that
$Y$ does not satisfy $\binom{\Omega}{\Tau}$.
Together with the remarks preceding Problem \ref{borhered}, we have that the
only classes (in addition to those in Problem \ref{borhered})
for which the hereditariness problem is not settled
are the following ones.

\begin{prob}[\Cite{Problem 4}{ideals}]\label{bortauhered}
Is any of the properties
$\sone(\BT,\BG)$, $\sone(\BT,\BT)$, $\sone(\BT,\BO)$, $\sone(\BT,\B)$,
$\sfin(\BT,\BT)$, or $\sfin(\BT,\BO)$, hereditary?
\end{prob}

Here are the open problems regarding unions.

\begin{prob}[\cite{AddQuad}]
Is any of the properties
$\sone(\Tau,\Tau)$, $\sfin(\Tau,\Tau)$,
$\sone(\Gamma,\Tau)$, $\sfin(\Gamma,\Tau)$,
and $\ufin(\Gamma,\Tau)$ (or any of their Borel
versions) preserved under taking finite unions?
\end{prob}

And here are the open problems regarding powers.

\begin{prob}\label{tauprod}
Is any of the properties
\be
\i $\sone(\Omega,\Tau)$, or $\sfin(\Omega,\Tau)$,
\i $\sone(\Tau,\Gamma)$, $\sone(\Tau,\Tau)$, $\sone(\Tau,\Omega)$,
$\sfin(\Tau,\Tau)$, or $\sfin(\Tau,\Omega)$,
\ee
preserved under taking finite powers?
\end{prob}
The answer to (1) is positive if it is for Problem \ref{omtau}.

\section{Groupable covers}\label{groupbl}

Groupability notions for covers appear naturally in the studies of selection principles
\cite{FunRez, coc7, coc8, hureslaloms} and have various notations.
Scheepers has standardized the notations in \cite{LecceSurvey}.
We use Scheepers' notation, but take a very minor variant of his definitions
which allows more simple definitions and does not make a difference in
any of the theorems proved in the literature (since we always consider covers which
do not contain finite subcovers).

\forget
An open cover $\mathcal{U}$ of a space $X$ is:
\be
\i \emph{$\gamma$-groupable}\footnote{In the cited papers this is simply called \emph{groupable}.}
if there is a partition of $\cU$ into finite sets,
$\cU = \Union_{n\in\N}\cF_n$, such  that for each $x\in X$ and all but finitely many $n$, $x\in \cup\cF_n$.
\i \emph{$\tau$-groupable}
if there is a partition of $\cU$ into finite sets, $\cU = \Union_{n\in\N}\cF_n$,
such  that for each $x,y\in X$, at least one of the sets
$\{n : x\in\cup\cF_n, y\nin\cup\cF_n\}$ or $\{n : y\in\cup\cF_n, x\nin\cup\cF_n\}$ is finite.
\i \emph{$\w$-groupable}\footnote{In \cite{coc8} this is called \emph{weakly groupable}.
Caution: The term ``$\w$-groupable'' is used in a different meaning in \cite{FunRez}.}
if there is a partition of $\cU$ into finite sets, $\cU = \Union_{n\in\N}\cF_n$,
such that for each finite subset $F$ of $X$, there exists $n$ such that $F\sbst\cup\cF_n$.
\ee
\emph{If $\cU$ does not contain a finite subcover} then:
\forgotten

Let $\xi$ be $\gamma$, $\tau$, or $\omega$.
A cover $\cU$ of $X$ is \emph{$\xi$-groupable} if
there is a partition of $\cU$ into finite sets, $\cU = \Union_{n\in\N}\cF_n$, such  that
$\seq{\cup\cF_n}$ is a $\xi$-cover of $X$.
Denote the collection of $\xi$-groupable open covers of $X$ by $\O^{\xi\gp}$.
Then $\O^{\gamma\gp}\sbst\O^{\tau\gp}\sbst\O^{\w\gp}$.
Observe that we must require in the definitions that the elements $\cF_n$ are disjoint,
as otherwise any cover of $X$ would be $\gamma$-groupable.

Recall that $\ufin(\Gamma,\O)\Iff\sfin(\Omega,\O)$.
In \cite{coc7} it is proved that $\ufin(\Gamma,\Gamma)\Iff\sfin(\Omega,\O^{\gamma\gp})$,
and in \cite{coc8} it is proved that $\ufin(\Gamma,\Omega)\Iff\sfin(\Omega,\O^{\w\gp})$.
\begin{prob}
Is $\ufin(\Gamma,\Tau)$ equivalent to $\sfin(\Omega,\O^{\tau\gp})$?
\end{prob}
In \cite{strongdiags} it was pointed out that $\sone(\Omega,\O^{\w\gp})$ is strictly stronger
than $\sone(\Omega,\O)$ (which is the same as $\sone(\O,\O)$).
The following problem remains open.
\begin{prob}[\Cite{Problem 4}{coc8}]\label{s1wgp}
Is $\sone(\Omega,\Omega)$ equivalent to $\sone(\Omega,\O^{\w\gp})$?
\end{prob}
If all powers of sets in $\sone(\Omega,\O^{\w\gp})$ satisfy $\sone(\O,\O)$, then
we get a positive answer to Problem \ref{s1wgp}.
In \cite{strongdiags} it is shown that
$\sone(\Omega,\O^{\w\gp})\Iff\ufin(\Gamma,\Omega)\cap\sone(\O,\O)$,
so the question can be stated in classical terms.
\begin{prob}
Is $\sone(\Omega,\Omega)$ equivalent to $\ufin(\Gamma,\Omega)\cap\sone(\O,\O)$?
\end{prob}

Surprisingly, it turns out that $\ufin(\Gamma,\Gamma)\Iff\binom{\Lambda}{\O^{\gamma\gp}}$ \cite{hureslaloms}.
\begin{prob}
\be
\i Is $\ufin(\Gamma,\Omega)$ equivalent to $\binom{\Lambda}{\O^{\w\gp}}$?
\i Is $\ufin(\Gamma,\Tau)$ equivalent to $\binom{\Lambda}{\O^{\tau\gp}}$?
\ee
\end{prob}
It is often the case that properties of the form $\Pi(\Omega,\fV)$ where
$\Omega\sbst\fV$ are equivalent to $\Pi(\Lambda,\fV)$ \cite{coc1, coc2, coc7, coc8, strongdiags}.
But we do not know the answer to the following simple question.
\begin{prob}
Is $\binom{\Omega}{\O^{\gamma\gp}}$ equivalent to $\binom{\Lambda}{\O^{\gamma\gp}}$?
\end{prob}

Let $\fU$ be a family of covers of $X$.
Following \cite{coc7}, we say that a cover $\cU$ of $X$ is \emph{$\fU$-groupable} if
there is a partition of $\cU$ into finite sets, $\cU = \Union_{n\in\N}\cF_n$, such  that
for each infinite subset $A$ of $\N$, $\{\cup\cF_n\}_{n\in A}\in\fU$.
Let $\fU^\gpbl$ be the family of $\fU$-groupable elements of $\fU$.
Observe that $\Lambda^\gpbl=\O^\gpbl=\O^{\gamma\gp}$.

In \cite{coc7} it is proved that $X$ satisfies $\sfin(\Omega,\Omega^{gp})$ if, and only if,
all finite powers of $X$ satisfy $\ufin(\Gamma,\Gamma)$, which we now know is the
same as $\binom{\Lambda}{\Lambda^{gp}}$.

\begin{prob}[\Cite{Problem 8}{reznicb}]\label{w/wgp}
Is $\sfin(\Omega,\Omega^{gp})$ equivalent to $\binom{\Omega}{\Omega^{gp}}$?
\end{prob}

In \cite{NSW, coc7} it is proved that $\sone(\Omega,\Lambda^\gpbl)\Iff\ufin(\Gamma,\Gamma)\cap\sone(\O,\O)$,
and is the same as the Gerlits-Nagy $(*)$ property.
Clearly, $\sone(\Omega,\Tau)$ implies $\sone(\Omega,\O)$, which is the same as $\sone(\O,\O)$.
Thus, a positive solution to Problem \ref{schtauprob} would imply
that the property $\sone(\Omega,\Tau)$ lies between the $\sone(\Omega,\Gamma)$ ($\gamma$-property)
and $\sone(\Omega,\Lambda^\gpbl)$ ($(*)$ property).

These notions of groupable covers are new and were not completely classified yet.
Some partial results appear in \cite{FunRez, coc7, coc8, hureslaloms}.
\begin{proj}
Classify the selection properties involving groupable covers.
\end{proj}
The studies of preservation of these properties under set theoretic operations
are also far from being complete. Some of the known results are quoted in \cite{LecceSurvey}.

\section{Splittability}

The following discussion is based on \cite{splittability}.
Assume that $\fU$ and $\fV$ are collections of covers of a space $X$.
The following property was introduced in \cite{coc1}.
\bi
\i[$\split(\fU,\fV)$:] Every cover $\cU\in\fU$ can be split
into two disjoint subcovers $\cV$ and $\cW$ which contain elements of $\fV$.
\ei
This property is useful in the Ramsey theory of thick covers.
Several results about these properties (where $\fU,\fV$ are collections of thick
covers) are scattered in the literature.
Some results relate these properties to classical
properties. For example, it is known that
the Hurewicz property and Rothberger's property each implies $\split(\Lambda,\Lambda)$, and that
the Sakai property (asserting that each finite power of $X$ has Rothberger's property)
implies $\split(\Omega,\Omega)$ \cite{coc1}.
It is also known that if all finite powers of $X$ have the Hurewicz property,
then $X$ satisfies $\split(\Omega,\Omega)$ \cite{coc7}.
Let $\CO$ denote the collection of all \emph{clopen} $\w$-covers of $X$.
By a recent characterization of the Reznichenko (or: weak Fr\'echet-Urysohn) property of $C_p(X)$
in terms of covering properties of $X$ \cite{SakaiRez},
the Reznichenko property for $C_p(X)$ implies that $X$ satisfies $\split(\CO,\CO)$.

\subsection{Classification}
If we consider this prototype with $\fU,\fV\in\{\Lambda,\Omega,\Tau,\Gamma\}$
we obtain the following $16$ properties.
\newcommand{\aru}{\ar[r]\ar[u]}
$$\xymatrix{
\split(\Lambda,\Lambda)\ar[r]&\split(\Omega,\Lambda)\ar[r]&\split(\Tau,\Lambda)\ar[r]&\split(\Gamma,\Lambda)\\
\split(\Lambda,\Omega)\aru&\split(\Omega,\Omega)\aru&\split(\Tau,\Omega)\aru&\split(\Gamma,\Omega)\ar[u]\\
\split(\Lambda,\Tau)\aru&\split(\Omega,\Tau)\aru&\split(\Tau,\Tau)\aru&\split(\Gamma,\Tau)\ar[u]\\
\split(\Lambda,\Gamma)\aru&\split(\Omega,\Gamma)\aru&\split(\Tau,\Gamma)\aru&\split(\Gamma,\Gamma)\ar[u]
}$$
But all properties in the last column are trivial in the sense that all sets
of reals satisfy them.
On the other hand, all properties but the top one in the first column imply $\binom{\Lambda}{\Omega}$ and are
therefore trivial in the sense that no infinite set of reals satisfies any of them.
Moreover, the properties $\split(\Tau , \Tau)$, $\split(\Tau , \Omega)$, and $\split(\Tau , \Lambda)$ are
equivalent.
It is also easy to see that $\split(\Omega,\Gamma) \Iff \binom{\Omega}{\Gamma}$, therefore
$\split(\Omega,\Gamma)$ implies $\split(\Lambda,\Lambda)$.
In \cite{splittability} it is proved that no implication can be added to the diagram in
the following problem, except perhaps the dotted ones.
\begin{prob}[\Cite{Problem 5.9}{splittability}]\label{splittabilityclassif}
Is the dotted implication (1) (and therefore (2) and (3)) in the following diagram true?
If not, then is the dotted implication (3) true?
$$
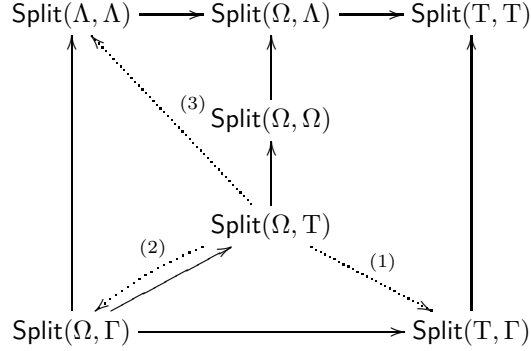
\xymatrix{
\split(\Lambda, \Lambda) \ar[r] & \split(\Omega, \Lambda) \ar[r] & \split(\Tau, \Tau)\\
                           & \split(\Omega, \Omega)\ar[u]\\
& \split(\Omega, \Tau)\ar[u]\ar@{.>}[dr]^{(1)}\ar@/_/@{.>}[dl]_{(2)}\ar@{.>}[uul]_{(3)}\\
\split(\Omega,\Gamma) \ar[uuu]\ar[ur]\ar[rr]     & & \split(\Tau,\Gamma)\ar[uuu]\\
}$$
\end{prob}
A positive answer to Problem \ref{omtau} would imply a positive answer to this problem.

\subsection{Preservation of properties}
We list briefly the only remaining problems concerning preservation of the splittability
properties mentioned in the last section under set theoretic operations.
All problems below are settled for the properties which do not appear in them.

\begin{prob}[\Cite{Problem 6.8}{splittability}]\label{SpLamLamAdd}
Is $\split(\Lambda,\Lambda)$ additive?
\end{prob}

\begin{prob}[\Cite{Problem 7.5}{splittability}]\label{BorelHered}
Is any of the properties $\split(\BO,\BL)$, $\split(\BO,\BO)$, $\split(\BT,\BT)$,
and $\split(\BT,\BG)$ hereditary?
\end{prob}

\begin{conj}[\Cite{Conjecture 8.7}{splittability}]\label{conj}
None of the classes $\split(\Tau,\Tau)$ and
$\binom{\Tau}{\Gamma}$ is provably closed under taking finite products.
\end{conj}

\begin{prob}[\Cite{Problem 8.8}{splittability}]\label{powcl}
Is any of the properties $\split(\Omega,\Omega)$, $\split(\Omega,\Tau)$, or
$\split(\Tau,\Tau)$ preserved under taking finite powers?
\end{prob}

\section{Ultrafilter-covers}

\subsection{The $\delta$-property}
The following problem is classical by now, but it is related to the
problems which follow.
For a sequence $\seq{X_n}$ of subsets of $X$, define
$\liminf X_n = \Union_m\bigcap_{n\ge m} X_n$.
For a family $\cF$ of subsets of $X$, $L(\cF)$ denotes
its closure under the operation $\liminf$.
The following definition appears
in the celebrated paper \cite{GN} just after that of the $\gamma$-property:
$X$ is a \emph{$\delta$-set} (or: has the \emph{$\delta$-property}) if
for each $\w$-cover $\cU$ of $X$, $X\in L(\cU)$.
Observe that if $\seq{U_n}$ is a $\gamma$-cover of $X$, then
$X=\liminf U_n$. Thus, the $\gamma$-property implies the $\delta$-property.
Surprisingly, the converse is still open.
\begin{prob}[\Cite{p.~160}{GN}]\label{GNdelta}
Is the $\delta$-property equivalent to the $\gamma$-property?
\end{prob}
The $\delta$-property implies Gerlits-Nagy $(*)$ property \cite{GN}, which is the same
as $\ufin(\Gamma,\Gamma)\cap\sone(\O,\O)$ (or $\sone(\Omega,\Lambda^\gpbl)$)
and implies $\sone(\Omega,\Omega)$ \cite{NSW}.
Miller (personal communication) suggested that if we could construct an increasing
sequence $\seq{X_n}$ of $\gamma$-sets whose union is not a $\gamma$-set, then the
union of these sets would be a $\delta$-set which is not a $\gamma$-set.

For a sequence $\seq{X_n}$ of subsets of $X$, define
$\plim X_n = \Union_{A\in p}\bigcap_{n\in A} X_n$.
For a family $\cF$ of subsets of $X$, $L_p(\cF)$ denotes
its closure under the operation $\plim$.
A space $X$ satisfies the $\delta_M$ property
if for each open $\w$-cover $\cU$ of $X$, there exists $p\in M$
such that $X\in S_p(\cU)$.
When $M=\{p\}$, we write $\delta_p$ instead of $\delta_{\{p\}}$.

The following problem is analogous to Problem \ref{GNdelta}.

\begin{prob}[\Cite{Problem 3.14}{GTpseq}]
Assume that $X$ satisfies $\delta_p$ for each ultrafilter $p$.
Must $X$ satisfy $\gamma_p$ for each ultrafilter $p$?
\end{prob}

\subsection{Sequential spaces}
A space $Y$ is \emph{sequential} if for each non-closed $A\sbst Y$
there exists $y\in Y\sm A$ and a sequence $\seq{a_n}$ in $A$
such that $\lim a_n = y$. This notion has a natural generalization.

An \emph{ultrafilter on $\N$} is a
family $p$ of subsets of $\N$ that
is closed under taking supersets, is closed under
finite intersections, does not contain the empty set as an element,
and for each $A\sbst\N$, either $A\in p$ or $\N\sm A\in p$.
An ultrafilter $p$ on $\N$ is \emph{nonprincipal} if it is not of the
form $\{A\sbst\N : n\in A\}$ for any $n$.
In the sequel, by \emph{ultrafilter} we always mean a \emph{nonprincipal ultrafilter on $\N$}.

For an ultrafilter $p$, $\O_p$ denotes the collection of open $\gamma_p$-covers of $X$,
that is, open covers $\cU$ that can be enumerated as $\seq{U_n}$ where
$\{n : x\in U_n\}\in p$ for all $x\in X$. The property $\binom{\Omega}{\O_p}$ is called
\emph{the $\gamma_p$-property} in \cite{GTpseq}.

\begin{prob}[\Cite{Question 2.4}{GJ}]
Is the property $\binom{\Omega}{\O_p}$ additive for each ultrafilter $p$?
\end{prob}

\begin{prob}[\Cite{Problem 3.14(2)}{GTpseq}]
Assume that $X$ satisfies $\binom{\Omega}{\O_p}$ for each ultrafilter $p$.
Must $X$ satisfy $\binom{\Omega}{\Gamma}$?
\end{prob}
In \cite[Theorem 3.13]{GTpseq} it is shown that the answer to this problem is positive
under an additional set theoretic hypothesis.

For an ultrafilter $p$, we write $x = \plim{} x_n$ when for each neighborhood $U$ of $x$,
$\{n : x_n\in U\}\in p$.
A space $Y$ is \emph{$p$-sequential} if we replace
$\lim$ by $\plim$ in the definition of sequential.

\begin{prob}[\cite{GT1}, \Cite{Problem 0.10}{GTpseq}]
Assume that $C_p(X)$ is $p$-sequential. Must $X$ satisfy $\binom{\Omega}{\O_p}$?
\end{prob}

Kombarov \cite{Km} introduced the following two generalizations of $p$-sequentiality:
Let $M$ be a collection of ultrafilters.
$Y$ is \emph{weakly $M$-sequential} if for
each non-closed $A\sbst Y$
there exists $y\in Y\sm A$ and a sequence $\seq{a_n}$ in $A$
such that $\plim a_n = y$ for \emph{some} $p\in M$.
$Y$ is \emph{strongly $M$-sequential} if \emph{some} is replaced
by \emph{for all} in the last definition.

\begin{prob}[\Cite{Problem 0.6 (reformulated)}{GTpseq}]
Assume that $X$ satisfies the $\delta_M$ property.
Must $C_p(X)$ be weakly $M$-sequential?
\end{prob}

\part{Applications}

\section{Infinite game Theory}

Each selection principle has a naturally associated game.
In the game $\gone(\fU,\fV)$ ONE chooses in the $n$th inning an element
$\cU_n$ of $\fU$ and then TWO responds by choosing $U_n\in \cU_n$.
They play an inning per natural number.
A play $(\cU_0, U_0, \cU_1, U_1\dots)$ is won by TWO if
$\seq{U_n}\in\fV$; otherwise ONE wins.
The game $\gfin(\fU,\fV)$ is played similarly, where TWO responds with
finite subsets $\cF_n\sbst\cU_n$ and wins if $\Union_{n\in\N}\cF_n\in\fV$.

Observe that if ONE does not have a winning strategy in
$\gone(\fU,\fV)$ (respectively, $\gfin(\fU,\fV)$), then $\sone(\fU,\fV)$
(respectively, $\sfin(\fU,\fV)$) holds.
The converse is not always true; when it is true,
the game is a powerful tool for studying the combinatorial
properties of $\fU$ and $\fV$ -- see, e.g., \cite{coc7},
\cite{coc8}, and references therein.

Let $\cD$ denote the collection of all families $\cU$ of open sets in $X$
such that $\cup\cU$ is dense in $X$.
In \cite{BerJu}, Berner and Juh\'asz introduce the open-point game,
which by \cite{coc6} is equivalent to $\gone(\cD,\cD)$ in the sense
that a player has a winning strategy in the open-point game on $X$
if, and only if, the \emph{other} player has a winning strategy in
$\gone(\cD,\cD)$.
\begin{prob}[\Cite{Question 4.2}{BerJu}, \Cite{footnote 1}{lengthdiags}]
Does there exist in ZFC a space in which $\gone(\cD,\cD)$ is undetermined?
\end{prob}

$\DO$ is the collection of all $\cU\in\cD$ such that for each $U\in\cU$,
$X\not\sbst\cl{U}$, and for each finite collection $\cF$ of open sets,
there exists $U\in\cU$ which intersects all members of $\cF$.

Problem \ref{sfindodo} is not a game theoretical one, but it is related to Problem
\ref{nextone} which is a game theoretic problem.
If all finite powers of $X$ satisfy $\sone(\cD,\cD)$ (respectively, $\sfin(\cD,\cD)$),
then $X$ satisfies $\sone(\DO,\DO)$ (respectively, $\sfin(\DO,\DO)$) \cite{coc5}.
If the other direction also holds, then the answer to the following is positive.
\begin{prob}[\cite{coc5}]\label{sfindodo}
Are the properties $\sone(\DO,\DO)$ or $\sfin(\DO,\DO)$ preserved under taking finite powers?
\end{prob}
The answer is ``Yes'' for a nontrivial family of spaces -- see \cite{coc5}.
A positive answer for this problem implies a positive answer to the following one.
If each finite power of $X$ satisfies $\sone(\cD,\cD)$ (respectively, $\sfin(\cD,\cD)$),
then ONE has no winning strategy in $\gone(\DO,\DO)$ (respectively, $\gfin(\DO,\DO)$) \cite{coc5}.
\begin{prob}[\Cite{Problem 3}{coc5}]\label{nextone}
Is any of the properties $\sone(\DO,\DO)$ or $\sfin(\DO,\DO)$
equivalent to ONE not having a winning strategy in the corresponding game?
\end{prob}

Let $\cK$ denote the families $\cU\in\cD$ such that $\{\cl{U} : U\in\cU\}$ is a cover of $X$.
In \cite{tkachuk} Tkachuk shows that \CH{} implies that ONE has a winning
strategy in $\gone(\cK,\cD)$ on any space of uncountable cellularity.
In \cite{coc4}, Scheepers defines $\fj$ as the minimal cardinal $\kappa$
such that ONE has a winning strategy in $\gone(\cK,\cD)$ on
each Tychonoff space with cellularity at least $\kappa$, and shows that
$\cov(\M)\le\fj\le\non(SMZ)$.
\begin{prob}[\Cite{Problem 1}{coc4}]
Is $\fj$ equal to any standard cardinal characteristic of the continuum?
\end{prob}
Scheepers conjectures that $\fj$ is not provably equal to $\cov(\M)$, and not to $\non(SMZ)$ either.

\subsection{Strong selection principles and games}

The following prototype of selection hypotheses is described in \cite{strongdiags}.
Assume that $\seq{\fU_n}$ is a sequence of collections of covers of a space $X$,
and that $\fV$ is a collection of covers of $X$.
Define the following selection hypothesis.
\begin{itemize}
\item[$\sone(\seq{\fU_n},\fV)$:]
For each sequence $\seq{\cU_n}$ where $\cU_n\in\fU_n$ for each $n$,
there is a sequence
$\seq{U_n}$ such that $U_n\in\cU_n$ for each $n$, and $\seq{U_n}\in\fV$.
\end{itemize}
Similarly, define $\sfin(\seq{\fU_n},\fV)$.
A cover $\cU$ of a space $X$ is an \emph{$n$-cover} if each
$n$-element subset of $X$ is contained in some member of $\cU$.
For each $n$ denote by $\O_n$ the collection of all
open $n$-covers of a space $X$.
Then $X$ is a \emph{strong $\gamma$-set} according to the definition of Galvin-Miller \cite{GM}
if, and only if, $X$ satisfies $\sone(\seq{\O_n},\Gamma)$ \cite{strongdiags}.
It is well known that the strong $\gamma$-property is strictly stronger than
the $\gamma$-property, and is therefore not equivalent to any of the classical properties.
However, for almost any other pair $(\seq{\fU_n},\fV)$, $\sone(\seq{\fU_n},\fV)$ and
$\sfin(\seq{\fU_n},\fV)$ turns out equivalent to some classical property \cite{strongdiags}.
The only remaining problem is the following.

\begin{conj}[\Cite{Conjecture 1}{strongdiags}]\label{strongOmTau}
$\sone(\seq{\O_n},\Tau)$ is strictly stronger than $\sone(\Omega,\Tau)$.
\end{conj}

If this conjecture is false, then we get a negative answer to Problem 13 of \cite{Hdim}.

As in the classical selection principles, there exist game-theoretical
counterparts of the new prototypes of selection principles \cite{strongdiags}.
Define the following games between two players, ONE and TWO,
which have an inning per natural number.
$\gone(\seq{\fU_n},\fV)$:
In the $n$th inning, ONE chooses an element $\cU_n\in\fU_n$, and
TWO responds with an element $U_n\in\cU_n$. TWO wins if
$\seq{U_n}\in\fV$; otherwise ONE wins.
$\gfin(\seq{\fU_n},\fV)$:
In the $n$th inning, ONE chooses an element $\cU_n\in\fU_n$, and
TWO responds with a finite subset $\cF_n$ of $\cU_n$. TWO wins if
$\Union_{n\in\N}\cF_n\in\fV$; otherwise ONE wins.

In \cite{strongdiags} it is proved that for $\fV\in\{\Lambda,\O^{\w\gp},\O^{\gamma\gp}\}$,
ONE does not have a winning strategy in $\gfin(\seq{\O_n},\fV)$ if, and only if,
$\sfin(\Omega,\fV)$ holds, and the analogous result is proved for $\gone$ and $\sone$.
In the case of $\gone$ and $\sone$, the assertion also holds for
$\fV\in\{\Omega,\Omega^{gp}\}$.

\begin{prob}
Assume that $\fV\in\{\Omega,\Omega^{gp}\}$.
Is it true that ONE does not have a winning strategy in
$\gfin(\seq{\O_n},\fV)$ if, and only if, $\sfin(\Omega,\fV)$ holds?
\end{prob}

The most interesting problem with regards to these games seems to be the following.
\begin{prob}[\Cite{Problem 5.16}{strongdiags}]
Is it true that $X$ is a strong $\gamma$-set (i.e., satisfies $\sone(\seq{\O_n},\Gamma)$)
if, and only if, ONE has no winning strategy in the game $\gone(\seq{\O_n},\Gamma)$?
\end{prob}
A positive answer would give the first game-theoretic characterization of
the strong $\gamma$-property.

\section{Ramsey Theory}

\subsection{Luzin sets}
Recall that $\cK$ is the collection of families $\cU$ of open sets such that $\{\cl{U}:U\in\cU\}$
is a cover of $X$. Let $\KO$ be the collection of all $\cU\in\cK$
such that no element of $\cU$ is dense in $X$, and for each
finite $F\sbst X$, there exists $U\in\cU$ such that
$F\sbst\cl{U}$.
In \cite{SchLuzin} it is proved that
$X$ satisfies $\KO\rightarrow(\cK)^2_2$, then $X$ is a Luzin set.
\begin{prob}[\Cite{Problem 4}{SchLuzin}]
Does the partition relation $\KO\rightarrow(\cK)^2_2$ characterize
Luzin sets?
\end{prob}

\subsection{Polarized partition relations}
The symbol
$$\(\mx{\fU_1\\\fU_2}\)\to\left[\mx{\fV_1\\\fV_2}\right]_{k/<\ell}$$
denotes the property that for each
$\cU_1\in\fU_1$, $\cU_2\in\fU_2$, and $k$-coloring
$f:\cU_1\x\cU_2\to\{1,\dots,k\}$ there are
$\cV_1\sbst\cU_1$, $\cV_2\sbst\cU_2$ such that $\cV_1\in\fV_1$ and $\cV_2\in\fV_2$,
and a set of less than $\ell$ colors $J$ such that
$f[\cV_1\x\cV_2]\sbst J$.

$\sone(\Omega,\Omega)$ implies
{\tiny $\(\mx{\Omega\\\Omega}\)\to\left[\mx{\Omega\\\Omega}\right]_{k/<3}$},
which in turn implies $\sfin(\Omega,\Omega)$ as well as $\split(\Omega,\Omega)$
(see Section \Ref{split} for the definition of the last property).
Consequently, the critical cardinality of this partition relation
lies between $\cov(\M)$ and $\min\{\d, \u\}$ \cite{PolarPar}.

\begin{prob}[\Cite{Problem 1}{PolarPar}]
Is {\tiny $\(\mx{\Omega\\\Omega}\)\to\left[\mx{\Omega\\\Omega}\right]_{k/<3}$}
equivalent to $\sone(\Omega,\Omega)$?
And if not, is its critical cardinality equal to that of $\sone(\Omega,\Omega)$
(namely, to $\cov(\M)$)?
\end{prob}

\section{Function spaces and Arkhangel'ski\v{i} duality theory}

The set of all real-valued functions on $X$, denoted
$\R^X$, is considered as a power of the real line and is endowed
with the Tychonoff product topology.
$C_p(X)$ is the subspace of $\R^X$ consisting of the continuous
real-valued functions on $X$. The topology of $C_p(X)$ is known as the
\emph{topology of pointwise convergence}.
The constant zero element of $C_p(X)$ is denoted $\mathbf{0}$.

\subsection{$s_1$ spaces and sequence selection properties}
In a manner similar to the observation made in Section 3 of \cite{wqn},
a positive solution to Problem \ref{addGGb} should imply a positive solution
to the following problem.
For subset $A\subset X$ we denote
\begin{eqnarray*}
s_0(A)&=&A,\quad s_{\xi}(A)=\{\lim_{n\to\infty} x_n:x_n\in
\bigcup_{\eta<\xi}s_{\eta}(A)\mbox{\ for each\ }n\in\N\},\\
\sigma(A)&=&\min\{\xi:s_{\xi}(A)=s_{\xi+1}(A)\}.
\end{eqnarray*}
Let $\Sigma(X)=\sup\{\sigma(A) : A\subseteq X\}$.
Fremlin \cite{FremlinSeq} proved that $\Sigma(C_p(X))$ must be $0$, $1$, or $\omega_1$.
If $\Sigma(C_p(X))=1$ then we say that $X$ is an \emph{$s_1$-space}.
\begin{prob}[{Fremlin \cite[Problem 15(c)]{FremlinSeq}}]\label{s1add}
Is the union of less than $\b$ many $s_1$-spaces an $s_1$-space?
\end{prob}

A sequence $\seq{f_n}\sbst C_p(X)$ converges \emph{quasi-normally} to a function $f$ on $X$
\cite{Ba} if there exists a sequence of positive reals $\seq{\epsilon_n}$
converging to $0$ such that for each $x\in X$
$|f_n(x)-f(x)|<\epsilon_n$ for all but finitely many $n$.
$X$ is a \emph{$\wQN$-space} \cite{BRR1} if
each sequence in $C_p(X)$ which converges to $\mathbf{0}$,
contains a quasi-normally convergent subsequence.

Finally, $C_p(X)$ has the \emph{sequence selection property (SSP)}
if for each sequence $\seq{\{f^n_k\}_{k\in\N}}$ of sequences in $C_p(X)$,
where each of them converges to $\mathbf{0}$,
there exists a sequence $\seq{k_n}$ such that
$\seq{f^n_{k_n}}$ converges to $\mathbf{0}$.
This is equivalent to Arkhangel'ski\v{i}'s $\alpha_2$ property of $C_p(X)$.

In \cite{wqn, FrwQN} it is shown that $s_1$ (for $X$), $\wQN$ (for $X$),
and SSP (for $C_p(X)$) are all equivalent.
This and other reasons lead to suspecting that all these equivalent properties
are equivalent to a standard selection hypothesis.
In \cite{wqn}, Scheepers shows that
$\sone(\Gamma,\Gamma)$ implies being an $\wQN$-space.
\begin{conj}[{Scheepers \cite[Conjecture 1]{wqn}}]
For sets of reals, $\wQN$ implies $\sone(\Gamma,\Gamma)$.
\end{conj}
If this conjecture is true, then Problems \ref{addGGb}
and \ref{s1add} coincide.

A space has \emph{countable fan tightness} if for each $x\in X$, if
$\seq{A_n}$ is a sequence of subsets of $X$ such that for each
$n$ $x\in\cl{A_n}$, then there are finite subsets $F_N\sbst A_n$, $n\in\N$,
such that $x\in\cl{\Union_n F_n}$. This property is due
to Arkhangel\'ski\v{i}, who proved in \cite{A} that
$C_p(X)$ has countable fan tightness if, and only if, every
finite power of $X$ satisfies $\ufin(\Gamma,\O)$
(this is the same as $\sfin(\Omega,\Omega)$).

The \emph{weak sequence selection property} for $C_p(X)$
\cite{CpHure} is defined as the SSP with the difference that
we only require that $\mathbf{0}\in\cl{\{f^n_{k_n} : n\in\N\}}$.

\begin{prob}[\Cite{Problem 1}{CpHure}]\label{cft}
Does countable fan tightness of $C_p(X)$ imply the weak sequence selection property?
\end{prob}

The \emph{monotonic sequence selection property} is defined like the SSP
with the additional assumption that for each $n$ the sequence $\{f^n_k\}_{k\in\N}$
converges pointwise \emph{monotonically} to $\mathbf{0}$.

\begin{prob}[\Cite{Problem 2}{CpHure}]\label{mss}
Does the monotonic sequence selection property of
$C_p(X)$ imply the weak sequence selection property?
\end{prob}

\section{The weak Fr\'echet-Urysohn property and Pytkeev spaces}

Recall that a topological space $Y$ has the Fr\'echet-Urysohn
property if for each subset $A$ of $Y$ and each
$y\in\cl{A}$, there exists a sequence $\seq{a_n}$
of elements of $A$ which converges to $y$.
If $y\nin A$ then we may assume that the elements
$a_n$, $n\in\N$, are distinct.
The following natural generalization of this property
was introduced by Reznichenko \cite{MaTi}:
$Y$ satisfies the \emph{weak Fr\'echet-Urysohn} property if
for each subset $A$ of $Y$ and each element
$y$ in $\cl{A}\sm A$, there exists a countably infinite
pairwise disjoint collection $\cF$ of finite subsets of $A$
such that for each neighborhood $U$ of $y$, $U\cap F\neq\emptyset$ for
all but finitely many $F\in\cF$.
In several works \cite{FunRez, coc7, SakaiRez}
this property appears as the \emph{Reznichenko} property.

In \cite{coc7} it is shown that
$C_p(X)$ has countable fan tightness as well as Reznichenko's property if, and only if,
each finite power of $X$ has the Hurewicz covering property.
Recently, Sakai found an exact dual of the Reznichenko property:
An open $\w$-cover $\cU$ of $X$ is
\emph{$\w$-shrinkable} if for each $U \in \cU$ there exists
a closed subset $C_U\sbst U$ such that
$\{C_U : U \in \cU\}$ is an $\w$-cover of $X$.
Then $C_p(X)$ has the Reznichenko property if, and only if, each $\w$-shrinkable open $\w$-cover of $X$
is $\w$-groupable \cite{SakaiRez}.
Thus if $X$ satisfies $\binom{\Omega}{\Omega^{gp}}$,
then $C_p(X)$ has the Reznichenko property. The other direction is not clear.

\begin{prob}[\Cite{Question 3.5}{SakaiRez}, \cite{reznicb}]
Is it true that $C_p(X)$ has the Reznichenko property if, and only if,
$X$ satisfies $\binom{\Omega}{\Omega^\gpbl}$?
\end{prob}

Another simply stated problem is the following.
\begin{prob}[\Cite{Question 3.6}{SakaiRez}]\label{sakaiprob}
Does $C_p(\NN)$ have the Reznicenko property?
\end{prob}

For a nonprincipal filter $\cF$ on $\N$ and a finite-to-one function
$f:\N\to\N$, $f(\cF):=\{A\sbst\N : f^{-1}[A]\in\cF\}$ is again a nonprincipal
filter on $\N$.
A filter $\cF$ on $\N$ is \emph{feeble} if
there exists a finite-to-one function $f$ such that
$f(\cF)$ consists of only the cofinite sets.
By Sakai's Theorem, if $C_p(X)$ has the Reznichenko property then $X$
satisfies $\binom{\CO}{\CO^\gpbl}$. In \cite{reznicb} it is shown that
$\binom{\CO}{\CO^\gpbl}$ is equivalent to the property that no continous
image of $X$ in the Rothberger space $\roth$ is a subbbase for a non-feeble
filter. Thus, if a subbase for a non-feeble filter cannot be a continuous image of $\NN$,
then the answer to Problem \ref{sakaiprob} is negative.

A family $\cP$ of subsets of of a space $Y$ is a \emph{$\pi$-network}
at $y\in Y$ if every neighborhood of $y$ contains some element of
$\cP$.
$Y$ is a \emph{Pytkeev space} if for each $y\in Y$ and
$A\sbst Y$ such that $y\in\cl{A}\sm A$,
there exists a countable $\pi$-network at $y$
which consists of infinite subsets of $A$.
In \cite{SakaiRez} it is proved that
$C_p(X)$ is a Pytkeev space if, and only if,
for each $\omega$-shrinkable open $\omega$-cover $\cU$ of $X$
there exist subfamilies $\cU_n\sbst\cU$, $n\in\N$, such that
$\bigcap_{n\in\N}\cU_n$ is an $\omega$-cover of $X$.
\begin{prob}[\Cite{Question 2.8}{SakaiRez}]
Can the term ``$\omega$-shrinkable'' be removed from Sakai's
characterization of the Pytkeev property of $C_p(X)$?
\end{prob}
If all finite powers of $X$ satisfy $\ufin(\Gamma,\Gamma)$,
then every open $\w$-cover of $X$ is $\w$-shrinkable \cite{SakaiRez},
thus a positive solution to the following problem would suffice.
\begin{prob}[\Cite{Question 2.9}{SakaiRez}]\label{pythu}
Assume that $C_p(X)$ is a Pytkeev space. Is it true that all finite powers of
$X$ satisfy $\ufin(\Gamma,\Gamma)$?
\end{prob}
Let $I=[0,1]$ be the closed unit interval in $\R$.
As all finite powers of $I$ are compact, and $C_p(I)$ is not a Pytkeev space \cite{SakaiRez},
the converse of Problem \ref{pythu} is false.

\textsc{Notes added in proof.}
Zdomskyy, in a series of recent works,
settled (or partially settled) some of the problems
mentioned in the paper, the answers being:
``Yes'' for Problem \ref{MnotH},
``No'' for Problem \ref{addM},
and ``Consistently yes'' for Problem \ref{SpLamLamAdd}.

Sakai showed that the answer for Problem \ref{w/wgp}
is ``No'', in the following strong sense: in his paper
\emph{Weak Fr\'echet-Urysohn property in function spaces},
it is proved that every analytic set of reals
(and, in particular, the Baire space $\NN$) satisfies $\binom{\BO}{\BO^\gpbl}$.
But we know that $\NN$ does not even satisfy Menger's property $\ufin(\O,\O)$.
This also answers Problem \ref{sakaiprob} in the affirmative.

Sakai also settled Problems \ref{cft} and \ref{mss} in the negative,
in his paper \emph{The sequence selection properties of $C_p(X)$},
Topology and its Applications \textbf{154}, 552--560.

The paper: H.\ Mildenberger, S.\ Shelah, and B.\ Tsaban,
\emph{The combinatorics of $\tau$-covers}
(~\arx{math.GN/0409068}~) contains new results simplifying
some problems.
Project \ref{cards} is almost completely settled ($4$ out of the
$6$ cardinals are found, the two remaining ones are equal but still
unknown). Consequently, $21$ out of the $76$ potential implications in
Project \ref{implications} are ruled out, consult this paper for
the updated list of problems in this project.

Finally, the preset author's paper \emph{Some new directions in infinite-combinatorial topology}
(in: \textbf{Set Theory}, eds.\ J.\ Bagaria and S.\ Todor\v{c}evic,
Trends in Mathematics, Birkhauser, 2006, 225--255.) contains a light introduction to the field and
several problems not appearing in the current survey.


\begin{thebibliography}{00}
\bibitem{A}
A.\ V.\ Arkhangel'ski\v{i},
\emph{Hurewicz spaces, analytic sets and fan tightness of function spaces},
Soviet Mathematical Doklady \textbf{33} (1986),
396--399.

\bibitem{coc8}
L.\ Babinkostova, Lj.\ D.\ R.\ Kocinac, and M.\ Scheepers,
\emph{Combinatorics of open covers (VIII)}, Topology and its
Applications \textbf{140} (2004),
15--32.

\bibitem{covM2}
T.\ Bartoszy\'nski and H.\ Judah,
\emph{On cofinality of the smallest covering of the real line by meager sets II},
Proceedings of the American Mathematical Society \textbf{123} (1995),
1879--1885.

\bibitem{jubar}
T.\ Bartoszy\'nski and H.\ Judah,
Set Theory: On the structure of the real line, A.\ K.\ Peters,
Massachusetts: 1995.

\bibitem{huremen2}
T.\ Bartoszy\'nski, S.\ Shelah, and B.\ Tsaban,
\emph{Additivity properties of topological diagonalizations},
The Journal of Symbolic Logic \textbf{68} (2003),
1254--1260.

\bibitem{ideals}
T.\ Bartoszy\'nski and B.\ Tsaban,
\emph{Hereditary topological diagonalizations and the Menger-Hurewicz Conjectures},
Proceedings of the American Mathematical Society \textbf{134} (2006),
605--615.

\bibitem{Bergelson}
Vitali Bergelson,
\emph{Ergodic Ramsey Theory -- an update},
in: \textbf{Ergodic Theory of $\Z^d$-actions} (ed.\ M.\ Pollicott and K.\ Schmidt),
London Math.\ Soc.\ Lecture Note Series \textbf{228} (1996),
1--61.

\bibitem{BerJu}
A.\ Berner and I.\ Juha\'sz,
\emph{Point-picking games and HFD's},
in: Models and Sets, Proceedings of the Logic Colloquium 1983,
Springer Verlag LNM \textbf{1103} (1984),
53--66.

\bibitem{BlassHBK}
A.\ R.\ Blass,
\emph{Combinatorial cardinal characteristics of the continuum},
in: \textbf{Handbook of Set Theory} (M.\ Foreman, A.\ Kanamori, and M.\ Magidor, eds.),
Kluwer Academic Publishers, Dordrecht, to appear.

\bibitem{JORG}
J.\ Brendle,
\emph{Generic constructions of small sets of reals},
Topology and it Applications \textbf{71} (1996),
125--147.

\bibitem{Ba}
Z.\ Bukovsk\'a,
\emph{Quasinormal convergence},
Math.\ Slovaca \textbf{41} (1991),
137--146.

\bibitem{LBopenprobs}
L.\ Bukovsk\'y,
\emph{Not distinguishing convergences: Open problems},
SPM Bulletin \textbf{3},
2--3.

\bibitem{BRR1}
L.\ Bukovsk\'y, I.\ Rec\l{}aw, and M.\ Repick\'y,
\emph{Spaces not distinguishing pointwise and quasinormal convergence of real functions},
Topology and its Applications \textbf{41} (1991),
25--40.

\bibitem{BuHa}
L.\ Bukovsk\'y and J.\ Hale\v s,
\emph{On Hurewicz properties},
Topology and its Applications \textbf{132} (2003), 71--79.

\bibitem{vD}
E.\ K.\ van~Douwen,
\emph{The  integers  and  topology},
in: \textbf{Handbook of Set Theoretic Topology} (eds.\ K.\ Kunen  and  J.\ Vaughan),
North-Holland, Amsterdam: 1984,
111--167.

\bibitem{FremlinSeq}
D.\ H.\ Fremlin,
\emph{Sequential convergence in $C_p(X)$},
Commentationes Mathematicae Universitatis Carolinae \textbf{35} (1994),
371--382.

\bibitem{FrwQN}
D.\ H.\ Fremlin,
\emph{SSP and WQN},
preprint.

\bibitem{GM}
F.\ Galvin and A.\ W.\ Miller,
\emph{$\gamma$-sets and other singular sets of real numbers},
Topology and it Applications \textbf{17} (1984),
145--155.


\renewcommand{\i}{i}
\bibitem{GJ}
S.\ Garcia-Ferreira and W.\ Just,
\emph{Some remarks on the $\gamma_p$ property},
Questions and Answers in General Topology \textbf{17} (1999),
1--8.

\bibitem{GT1}
S.\ Garcia-Ferreira and A.\ Tamariz-Mascar\'ua,
\emph{$p$-Fr\'echet-Urysohn properties of function spaces},
Topology and its Applications \textbf{58} (1994),
157--172.

\bibitem{GTpseq}
S.\ Garcia-Ferreira and A.\ Tamariz-Mascar\'ua,
\emph{$p$-sequential like properties in function spaces},
Commentationes Mathematicae Universitatis Carolinae \textbf{35} (1994),
753--771.

\bibitem{GN}
J.\ Gerlits and Zs.\ Nagy,
\emph{Some properties of $C(X)$, I},
Topology and its Applications \textbf{14} (1982),
151--161.

\bibitem{HURE27}
W.\ Hurewicz,
\emph{\"Uber Folgen stetiger Funktionen},
Fundamenta Mathematicae \textbf{9} (1927),
193--204.

\bibitem{coc2}
W.\ Just, A.\ W.\ Miller, M.\ Scheepers, and P.\ J.\ Szeptycki,
\emph{The combinatorics of open covers II},
Topology and its Applications \textbf{73} (1996),
241--266.

\bibitem{Km}
A.\ P.\ Kombarov,
\emph{On a theorem of A.\ H.\ Stone},
Soveit Math.\ Dokl.\ \textbf{27} (1983),
544--547.

\bibitem{FunRez}
Lj.\ D.\ R.\ Ko\v{c}inac and M.\ Scheepers,
\emph{Function spaces and a property of Reznichenko},
Topology and it Applications \textbf{123} (2002),
135--143.

\bibitem{coc7}
Lj.\ D.\ R.\ Ko\v{c}inac and M.\ Scheepers,
\emph{Combinatorics of open covers (VII): Groupability},
Fundamenta Mathematicae \textbf{179} (2003),
131--155.

\bibitem{MaTi}
V.\ I.\ Malykhin and G.\ Tironi,
\emph{Weakly Fr\'echet-Urysohn and Pytkeev spaces},
Topology and its applications \textbf{104} (2000),
181--190.

\bibitem{Mil79Len}
A.\ W.\ Miller,
\emph{On the length of Borel hierarchies},
Annals of Mathematical Logic \textbf{16} (1979),
233--267.

\bibitem{MilSpec}
A.\ W.\ Miller,
\emph{Special subsets of the real line},
in: Handbook of Set-Theoretic Topology
(K.\ Kunen and J.\ E.\ Vaughan, eds.),
North Holland, Amsterdam: 1984,
201--233.

\bibitem{MilNonGamma}
A.\ W.\ Miller,
\emph{A Nonhereditary Borel-cover $\gamma$-set},
Real Analysis Exchange \textbf{29} (2003/4),
601--606.

\bibitem{NSW}
A.\ Nowik, M.\ Scheepers, and T.\ Weiss,
\emph{The algebraic sum of sets of real numbers with strong measure zero sets},
J.\ Symbolic Logic \textbf{63} (1998),
301--324.

\bibitem{pawlikowskireclaw}
J.\ Pawlikowski and I.\ Rec{\l}aw,
\emph{Parametrized Cicho\'n's diagram and small sets},
Fundamenta Mathematicae \textbf{147} (1995),
135--155.

\bibitem{luzinundetermined}
I.\ Rec\l{}aw,
\emph{Every Luzin set is undetermined in the point-open game},
Fundamenta Mathematicae \textbf{144} (1994),
43--54.

\bibitem{sakai}
M.\ Sakai,
\emph{Property $C''$ and function spaces},
Proceedings of the American Mathematical Society \textbf{104} (1988),
917--919.

\bibitem{SakaiRez}
M.\ Sakai,
\emph{The Pytkeev property and the Reznichenko property in function spaces},
Note di Matematica \textbf{22} (2003),
43--52.

\bibitem{coc1}
M.\ Scheepers,
\emph{Combinatorics of open covers I: Ramsey theory},
Topology and its Applications \textbf{69} (1996),
31--62.

\bibitem{CpHure}
M.\ Scheepers,
\emph{A sequential property of $ \mathsf{C}_p(X)$ and a covering property of Hurewicz},
Proceedings of the American Mathematical Society \textbf{125} (1997),
2789--2795.

\bibitem{coc4}
M.\ Scheepers,
\emph{Combinatorics of open covers (IV): subspaces of the Alexandroff double of the unit interval},
Topology and its Applications \textbf{83} (1998),
63--75.

\bibitem{SchLuzin}
M.\ Scheepers,
\emph{Lusin Sets},
Proceedings of the American Mathematical Society \textbf{127} (1999),
251--257.

\bibitem{wqn}
M.\ Scheepers,
\emph{Sequential convergence in $\mathsf{C}_p(X)$ and a covering property},
East-West Journal of Mathematics \textbf{1} (1999),
207--214.

\bibitem{lengthdiags}
M.\ Scheepers,
\emph{The length of some diagonalization games},
Arch.\ Math.\ Logic \textbf{38} (1999),
103--122.

\bibitem{coc5}
M.\ Scheepers,
\emph{Combinatorics of open covers (V): Pixley-Roy spaces and sets of reals, and omega-covers},
Topology and its Applications \textbf{102} (2000),
13--31.

\bibitem{coc6}
M.\ Scheepers,
\emph{Combinatorics of open covers (VI): Selectors for sequences of dense sets},
Quaestiones Mathematicae \textbf{22} (1999),
109--130.

\bibitem{alpha_i}
M.\ Scheepers,
\emph{$C_p(X)$ and Arkhangel'ski\v{i}'s $\alpha_i$ spaces},
Topology and its Applications \textbf{89} (1998),
265--275.

\bibitem{PolarPar}
M.\ Scheepers,
\emph{Open covers and the square bracket partition relation},
preprint.

\bibitem{CBC}
M.\ Scheepers and B.\ Tsaban,
\emph{The combinatorics of Borel covers},
Topology and its Applications \textbf{121} (2002),
357--382.
\arx{math.GN/0302322}

\bibitem{LecceSurvey}
M.\ Scheepers,
\emph{Selection principles and covering properties in topology},
Note di Matematica \textbf{22} (2003),
3--41.

\bibitem{ShTb768}
S.\ Shelah and B.\ Tsaban,
\emph{Critical cardinalities and additivity properties of combinatorial notions of smallness},
Journal of Applied Analysis \textbf{9} (2003),
149--162.
\arx{math.LO/0304019}

\bibitem{tkachuk}
V.\ V.\ Tkachuk,
\emph{Some new versions of an old game},
Commentationes Mathematicae Universitatis Carolinae \textbf{36} (1995),
177--196.

\bibitem{tau}
B.\ Tsaban,
\emph{A topological interpretation of $\mathfrak{t}$},
Real Analysis Exchange \textbf{25} (1999/2000),
391--404.
\arx{math.LO/9705209}

\bibitem{spm}
B.\ Tsaban (ed.),\\
SPM Bulletin \textbf{1} (2003), \arx{math.GN/0301011}\\
SPM Bulletin \textbf{2} (2003), \arx{math.GN/0302062}\\
SPM Bulletin \textbf{3} (2003), \arx{math.GN/0303057}\\
SPM Bulletin \textbf{4} (2003), \arx{math.GN/0304087}\\
SPM Bulletin \textbf{5} (2003), \arx{math.GN/0305367}\\
SPM Bulletin \textbf{6} (2003), \arx{math.GN/0312140}\\
SPM Bulletin \textbf{7} (2004), \arx{math.GN/0401155}\\
SPM Bulletin \textbf{8} (2004), \arx{math.GN/0403369}\\
SPM Bulletin \textbf{9} (2004), \arx{math.GN/0406411}\\
SPM Bulletin \textbf{10} (2004), \arx{math.GN/0409072}

\bibitem{splittability}
B.\ Tsaban,
\emph{The combinatorics of splittability},
Annals of Pure and Applied Logic \textbf{129} (2004),
107--130.
\\ \arx{math.LO/0212312}

\bibitem{tautau}
B.\ Tsaban,
\emph{Selection principles and the minimal tower problem},
Note di Matematica \textbf{22} (2003),
53--81.
\\ \arx{math.LO/0105045}

\bibitem{reznicb}
B.\ Tsaban,
\emph{The minimal cardinality where the Reznichenko property fails},
Israel Journal of Mathematics \textbf{140} (2004),
367--374.
\\ \arx{math.GN/0304024}

\bibitem{Hdim}
T.\ Weiss and B.\ Tsaban,
\emph{Topological diagonalizations and Hausdorff dimension},
Note di Matematica \textbf{22} (2003),
83--92.
\\ \arx{math.LO/0212009}

\bibitem{o-bdd}
B.\ Tsaban,
\emph{$o$-bounded groups and other topological groups with strong combinatorial properties},
Proceedings of the American Mathematical Society \textbf{134} (2006),
881--891.

\bibitem{strongdiags}
B.\ Tsaban,
\emph{Strong $\gamma$-sets and other singular spaces},
Topology and its Applications \textbf{153} (2005),
620--639.
\\ \arx{math.LO/0208057}

\bibitem{hureslaloms}
B.\ Tsaban,
\emph{The Hurewicz covering property and slaloms in the Baire space},
Fundamenta Mathematicae \textbf{181} (2004),
273--280.
\\ \arx{math.GN/0301085}

\Bc{AddQuad}{B. Tsaban}{Additivity numbers of covering properties}{Selection Principles and Covering Properties in Topology}
{L. Ko\v{c}inac, ed.}{Quaderni di Matematica 18, Seconda Universita di Napoli, Caserta}{2006}{245}{282}

\bibitem{prods}
B.\ Tsaban and T.\ Weiss,
\emph{Products of special sets of real numbers},
Real Analysis Exchange \textbf{30} (2004/5),
819--836.

\end{thebibliography}
\end{document}